\documentclass[12pt]{article}
\usepackage{amssymb,amsthm,amscd,amsmath,graphicx,color}

\newcommand{\R}{{\mathbb R}}

\newtheorem{theorem}{Theorem}[section]
\newtheorem{corollary}[theorem]{Corollary}
\newtheorem{lemma}[theorem]{Lemma}
\newtheorem{proposition}[theorem]{Proposition}
\newtheorem{definition}[theorem]{Definition}
\newtheorem{remark}[theorem]{Remark}

\begin{document}

\thispagestyle{empty}
\begin{flushright} \rm J. Geo. Phys. {\bf
  56/9} (2006), 1688-1708\end{flushright}
\par\bigskip\par
\vfill
\begin{center}
{\bfseries\Large
Two-dimensional Lagrangian singularities and bifurcations of gradient lines I}
\par\addvspace{20pt}
{\sc G. Marelli}
\par\medskip
Department of Mathematics, Kyoto University,\\
Kitashirakawa, Sakyo-ku, Kyoto 606-8502, Japan\end{center}
\vfill
\begin{quote} \footnotesize {\sc Abstract.}
Motivated by mirror symmetry,
we consider a Lagrangian fibration $X\rightarrow B$ and
Lagrangian maps $f:L\hookrightarrow X\rightarrow B$, 
when $L$ has dimension 2, exhibiting an unstable singularity, and study 
how their caustic changes, in a neighbourhood of the unstable singularity,
when slightly perturbed. 
The integral curves of $\nabla f_x$, for $x\in B$,
where $f_x(y)=f(y)-x\cdot y$, called ``gradient lines'',
are then introduced, and a study of them,
in order to analyse their bifurcation locus,
is carried out.
\end{quote}
\vfill
\leftline{\hbox to8cm{\hrulefill}}\par
{\footnotesize
\noindent Supported by a JSPS Postdoctoral fellowship
\par
 \par
\noindent\emph{E-Mail address:} {\tt marelli@kusm.kyoto-u.ac.jp}
}
\eject\thispagestyle{empty}\mbox{\ \ \ }\newpage\setcounter{page}{1}

\section{Introduction}

This is first of two papers motivated by the attempt of 
understanding some aspects
of the homological mirror conjecture, when we assume the existence of
dual torus fibrations. In it, we are concerned with the torus
fibration $T^{2n}\rightarrow T^n$, 
as a first step in the direction towards generic Lagrangian torus
fibrations. 
In this case mirror symmetry has
been studied, under certain hypothesis, in papers such as \cite{F2},
\cite{AP}, \cite{LYZ}, \cite{BMP1} and \cite{BMP2}, where the
idea that mirror symmetry is a kind of Fourier-Mukai transform has
been developed: given
a Lagrangian submanifold $L$ of $X$ supporting a local system, under certain
assumptions, a holomorphic bundle is obtained on a submanifold of
$X^\vee$.
In all these papers, a crucial hypothesis is that the caustic of
$L$ is empty, that is, the composition $L\hookrightarrow X\rightarrow
B$ has no critical points. This paper, instead, 
takes the first steps in the direction of including the
caustic. If $K\subset B$ denotes the caustic of $L$, we may think to
restrict the fibration to $B\setminus K$: now $L$ has no caustic
and we may apply what is known in this case and obtain a
holomorphic bundle on a certain submanifold of $X^\vee$ fibred over
$B\setminus K$; however we realize that the holomorphic structure
presents a monodromy which prevents from extending the
holomorphic bundle over the points of the caustic $K$.  As foreseen in
\cite{F2}, some quantum corrections must be performed
in order to extend the
holomorphic structure over points of the caustic.  Quantum
corrections or instanton effects are provided by pseudoholomorphic
discs in $X$ which bound $L$. Following \cite{F1}, the fibre over
$x\in B$ of the holomorphic
bundle on $X^\vee$ is constructed as Lagrangian intersection Floer
homology of $L$ and of the Lagrangian fibre of $X$ over $x$. This
approach is
equivalent to the Fourier-Mukai one when the caustic is empty,
but, unlike this,
has the advantage of naturally including pseudoholomorphic
discs. Assuming that, near $K$, 
Lagrangian intersection Floer homology is equivalent to
Morse homology defined through the generating function of $L$,
assumption which still must be
clarified and proved, enables us to study gradient lines of
$\nabla f_x$ instead of pseudo-holomorphic discs. 
This is the idea which leads the
development of this paper.
The theory of Lagrangian maps provides a classifications of Lagrangian
singularities: in dimension 2 only folds and cusps are generic and
stable; in dimension 3 other singularities appear, and so on. This
suggests us to start by considering the case when $L$ has dimension 2
and so $X$ is the torus fibration
$T^4\rightarrow T^2$. If $f$ is a (local) generating function of $L$,
we plan to study, in a neighbourhood of a point $x\in K$, the
gradient lines of the vector field $\nabla f_x$, where
$f_x(y)=f(y)-x\cdot y$, and their bifurcations, and with these to
construct Morse homology. Troubles are given by those singularities
which appear in dimension 2 as unstable, such as the elliptic
umbilic. For these we study what happens to the caustic and, in
the case of the elliptic umbilic, to the
bifurcation locus of gradient lines, when a small perturbation is added
to the
generating function $f$. In fact, in the Fukaya category, Lagrangian
maps are considered up to Hamiltonian equivalence, so we expect to
recover the case of $L$ having an unstable singularity by studying the
case of a Lagrangian submanifold $L'$ exhibiting a
stable singularity and Hamiltonian equivalent to $L$. The analysis of
possible phase portraits of $\nabla f_x$ should allow to construct the Morse
complex, while the study of reciprocal positions of the caustic and
bifurcation locus, providing morphisms of the Morse complex, should allow to
construct a bundle whose holomorphic strucure can be extended to the
caustic.

In this paper, after reviewing in section \ref{lagrmaps} 
some aspects about the 
classification
of Lagrangian singularities, we study in section
\ref{unstsing} how unstable critical
points of a Lagrangian map 
split when its generating function $f$ is slightly perturbed. 
We first consider a map whose
caustic is reduced to an elliptic umbilic, in a sense which we will
specify, and see that a small
perturbation modifies the caustic in a well-known curve known as tricuspoid. A
similar analysis is sketched for maps exhibiting other unstable 
singularities, such
as the hyperbolic umbilic, the swallow-tail and the parabolic
umbilic, deducing some ideas about the way the problem could be
faced. 
However, we recognize that, if we are interested in
application to Mirror Symmetry, in dimension 2 the relevant
singularities are the fold and the cusp, which are also stable, and the
elliptic umbilic, the hyperbolic umbilic and the swallow-tail, which, as said,
are stable and generic in dimension 3 though unstable in dimension 2. 
Singularities such as the
parabolic 
umbilic, stable and generic in dimension at least 4, 
become relevant only when studying the problem 
in dimension 3. The analysis of the gradient
lines of $\nabla f_x$ and their bifurcations occupy the whole
section \ref{bifgradlin}. 
This is essential to construct the Morse complex, 
in view of applications to Mirror Symmetry.
In particular we study which kind of bifurcations can occur 
in a family of vector 
fields exhibiting only saddles and nodes, 
and, given a bifurcation diagram, when there exists a family of
gradient vector fields providing that diagram.
In \cite{M}, the sequel to this paper, we will
analyse, in some specific cases, how the
bifurcation locus changes when the generating function $f$ is slightly
perturbed. 

\par\smallskip
{\bf Acknowledgements.} I wish to thank K.~Fukaya, whose suggestions
and 
help were decisive for
the achievement of all the results here expounded. 

I am thankful to JSPS (Japan Society for the 
Promotion of Science) which awarded me with a JSPS postdoctoral
fellowship at 
Kyoto University,
where this paper was written.

\section{Lagrangian submanifolds and their singularities}
\label{lagrmaps}

We recall some facts about Lagrangian 
submanifolds and 
their singularities, referring to
\cite{V} or \cite{AGZV} for details.

\subsection{Lagrangian maps}

Let $(X,\omega)$ be a symplectic $2n$-manifold, which will be
denoted simply by $X$, and $L$ an $n$-submanifold of $X$.

\begin{definition}
An immersion $g:L\rightarrow X$ is called 
Lagrangian immersion if $g^*\omega=0$. If $L\subset X$ and the
identical embedding is a Lagrangian immersion, then $L$ is called Lagrangian
submanifold.
\end{definition}

\begin{definition}
A Lagrangian bundle is a symplectic manifold $X$ endowed with a
structure of smooth
locally trivial bundle $\pi:X\rightarrow B$ over a base manifold $B$, 
all of whose fibres 
are Lagrangian submanifolds of $X$. 
\end{definition}

By  
Darboux's theorem,
any point of $X$ admits a neighbourhood with canonical
coordinates, that is, coordinates $(y_1,x_1,...,y_n,x_n)$ which are
both canonical symplectic coordinates of $X$ and such that the functions $x_i$
are constant along the fibres of the bundles. 

\begin{definition} 
Let $\pi:X\rightarrow B$ be a Lagrangian bundle and 
$g:L\rightarrow X$ a Lagrangian immersion. We call  Lagrangian map
the composite 
map $\pi\circ g:L
\rightarrow B$. The set $K$ of critical values of
$\pi\circ g$ is called caustic of the Lagrangian immersion $g$ (or of
the Lagrangian map $\pi\circ g$).
\end{definition}

If non-empty,
the caustic of a general Lagrangian map is an $(n-1)$-submanifold of $B$
with singularities. A classification of singularities of
Lagrangian maps is 
available and is obtained from the classification of 
singularities of smooth maps. Lagrangian maps are traced back to
smooth functions by means of their 
generating function: let $L\hookrightarrow X$ be 
a Lagrangian submanifold,
$p\in L$ a point, $\{y_i,x_i\}$ a system of canonical coordinates near $p$,
then there is a set of indices $J\subset\{1,...,n\}$ such that, if 
$I=\{1,...,n\}\setminus J$, then $\{y_i,x_j\}$ 
are local coordinates of $L$ near $p$, with $i\in I$ and $j\in J$, and
there exists a
smooth function $f$ in the variables $\{y_i,x_j\}$, defined up to addition
of a constant, such that $L$ is determined by the
equations
\begin{equation}
\label{genfun}
\left\{ \begin{aligned}
x_j & =\frac{\partial f}{\partial y_i}\\
y_i & =-\frac{\partial f}{\partial x_j}
\end{aligned}
\right.
\end{equation}
conversely, given a function $f$ as before, equations (\ref{genfun})
define a Lagrangian submanifold.
The function $f$ 
is called generating function of $L$.

\begin{definition}
Two Lagrangian bundles are said to be Lagrangian equivalent if there exists a
bundle diffeomorphism between them, taking fibres to fibres, and mapping one 
symplectic form to the other.
Analogously, two Lagrangian maps are said to be Lagrangian equivalent if there
exists a Lagrangian equivalence of the corresponding fibre bundles sending the
domain of the first map to that of the second.
\end{definition}

If two maps are Lagrangian equivalent then their caustic are
diffeomorphic. The converse of this statement is false. 

In studying the classification of Lagrangian singularities, it is more
convenient to enlarge the number of variables and describe a Lagrangian
germ by a function of the enlarged
set of variables and called generating family:
for a given Lagrangian germ, a generating families is not uniquely
determined, however
the class defining equivalent Lagrangian germs can be described.
If
$f(y_i,x_j)$ is the generating function of a germ of a Lagrangian submanifold
$L$, then 
$$F(z,x)=f(z_i,x_j)+<z_j,x_j>$$ 
is a generating family of $L$. Given
$F$, then
$L$ can be described as
the set
$$L=\{(y,x)~|~\exists z\quad with\quad 
\partial F/\partial z=0,\quad y=\partial F/\partial x\}$$
and its caustic $K$ as
$$K=\{x~|~\exists z\quad with\quad\partial F/\partial z=0,\quad 
det(\partial^2 F/\partial^2 x)=0\}$$

Let $D_0$ be the group of germs at $0$ of diffeomorphisms of $\R^n$ preserving 
$0$.

\begin{definition}
Two germs $f_1$ and $f_2$ of functions at $0$ are $D_0$-equivalent if there 
exists a germ $\phi\in D_0$ such that $f_1=f_2\circ\phi$.\\
Two germs $f_1$ and $f_2$ of functions at $0$ are stably $D_0$-equivalent if 
there  exists a germ $\phi\in D_0$ and a non-degenerate quadratic form $Q$ in 
additional variables such that $f_1=f_2\circ\phi+Q$.\\
\end{definition}

\begin{theorem}
\label{equiv}
Germs of Lagrangian maps are Lagrangian 
equivalent
if and only if their generating families are stably equivalent.
\end{theorem}

\subsection{The Whitney topology}
Let $X$ be a $2n$-symplectic manifold and $X\rightarrow B$ a
Lagrangian bundle.
Generating functions of Lagrangian maps are elements of $C^\infty(\R^n)$.
We endow the space of smooth function $C^\infty(\R^n)$
with
the Whitney $C^\infty$ topology (see also \cite{GG} and \cite{DM}):

\begin{definition}
For every non-negative integer $k$, and for every subset 
$U\subset J^k(\R^n)$, where
$J^k(\R^n)$ denotes the space of $k$-jets of smooth functions, let\\ 
$M(U)=\{f\in C^\infty(\R^n)~|~j^kf(\R^n)\subset U\}$. 
The family of sets $\{M(U)\}$
forms a basis for the Whitney $C^k$ topology on $C^\infty(\R^n)$. The
Whitney $C^\infty$ topology is the topology with basis
$W=\cup_{k=0}^\infty W_k$, where $W_k$ is the set of open subsets of
$C^\infty(\R^n)$ in the Whitney $C^k$ topology.
\end{definition}
 
Endowed with the Whitney $C^\infty$ topology, $C^\infty(\R^n)$ is a
Baire space, so every residual subset is dense. 

\begin{definition} 
A Lagrangian map is said to be Lagrangian stable if every nearby
Lagrangian map, in the Whitney topology, is Lagrangian equivalent to
it. 
\end{definition}

It can be proved that a germ of a Lagrangian map given by a generating
family $F$ is Lagrangian stable if and
only if $F$ is a versal deformation of $f(y_i,0)$, and that
its caustic is a component of the bifurcation 
set of its generating family.

\begin{definition}
A property $P$ of smooth functions in $C^\infty(\R^n)$ is generic 
if:\\ 
1. $C_P=\{f\in C^\infty(\R^n)~|~f~\textrm{satisfies}~P\}$ contains a
residual subset of $C^\infty(\R^n)$;\\ 
2. let $f\in C_P$ and suppose $g$ is
Lagrangian equivalent to $f$, then $g\in C_P$.
\end{definition}

A quasi-norm, and so a metric, generating the Whitney $C^\infty$
topology, can be defined on $C^\infty(\R^n)$
(see again \cite{GG} and
\cite{DM} for details), so that it makes sense to talk of  
small perturbations of a function $f\in C^\infty(\R^n)$.

\subsection{Classifications of Lagrangian singularities}
According to theorem \ref{equiv},
the problem of classifying Lagrangian singularities is reduced to
classify singularities of functions up to 
stably 
$D_0$-equivalence. The next theorem explains what happens in low
dimensions. For a list of normal forms see \cite{AGZV} or \cite{V}.

\begin{theorem}
The germs of generic Lagrangian maps $L\hookrightarrow X\rightarrow B$, with
$L$ of dimension $n\leq 5$, are stable and belong to a finite number of classes
of Lagrangian equivalence. When $n>5$ moduli appear, which in higher dimensions
become functional moduli. A classification of generic Lagrangian 
singularities exists for $n\leq 10$.
\end{theorem}

When $n\leq3$, the 
possible generating 
functions, denoted by letters $A$ or $D$, together with an index which 
represents the Milnor number, are:\\\\
$n\geq 1$
\begin{equation}
\label{gffold}
the~fold\quad A_2:~~~f(y_1)=y_1^3\quad
\end{equation}\\
$n\geq 2$
\begin{equation}
the~cusp\quad A_3:~~~f(y_1,x_2)=\pm y_1^4+x_2 y_1^2
\end{equation}
$n\geq 3$
\begin{equation}
the~swallow~tail\quad A_4:~~~f(y_1,x_2,x_3)=y_1^5+x_2 y_1^3+x_3 y_1^2
\end{equation}
\begin{equation}
\label{hugenfun}
the~hyperbolic~umbilic~or~purse\quad
D_4^+:~~~f(y_1,y_2,x_3)=y_1^3+y_1 y_2^2+x_3 y_1^2
\end{equation}
\begin{equation}
\label{eugenfun}
the~elliptic~umbilic~or~pyramid\quad
D_4^-:~~~f(y_1,y_2,x_3)=y_1^3-y_1 y_2^2+x_3 y_1^2
\end{equation}\\

\section{Perturbations of 2-dimensional unstable singularities}
\label{unstsing}
Let $X$ be a $4$-symplectic manifold and $X\rightarrow B$ a
Lagrangian bundle.
When Lagrangian submanifolds have dimension $2$,
only folds and cusps can appear locally as singularities of
generic stable Lagrangian maps, however other singularities can appear as non 
generic ones. In this case, such singularities are not stable and 
break in folds and cusps as a consequence of any generic perturbations.

Suppose that a Lagrangian map has an unstable critical point
at $p$ and we want to study how this singularity decomposes after a small
perturbation. 
To this purpose, 
consider
a small perturbation $f'$, which we can suppose supported on
a disc $D$ containing $p$.
This defines a new generating function $\tilde{f}=f+f'$ 
and a Lagrangian submanifold $\tilde{L}$; 
$L$ and $\tilde{L}$ coincide outside a compact subset 
$D'$ of $X$ and 
their caustics differ only in $f(D)\subset B$.

Being interested, as a first step, in a Lagrangian torus fibration with
2-dimensional smooth fibres, and since the decomposition of an unstable
singularity is a local problem, we can consider the Lagrangian
fibration $\R^4\rightarrow\R^2$. We use cooordinates $(x_1,x_2)$ on
the base and $(y_1,y_2)$ on the fibres.

\subsection{The elliptic umbilic}
\label{sseu}

We refer to
the generating function
\begin{equation}
\label{eugenfun2}
f(y_1,y_2)=\frac{1}{3}y_1^3-2y_1y_2^2
\end{equation}
defining the Lagrangian map
\begin{equation}
\label{elumbi}
\left\{ \begin{aligned}
x_1 & =y_1^2-y_2^2\\
x_2 & =-2y_1y_2
\end{aligned}
\right.
\end{equation}
as the elliptic umbilic in dimension 2.
It has an unique critical point, the origin $(0,0)$ of the
$(y_1,y_2)$-plane: it is neither a fold nor a cusp, so it is
unstable. The caustic is the subset $\{(0,0)\}$  
of the $(x_1,x_2)$-plane. To study how it splits when $f$ is slightly
perturbed, we add a perturbation $f'$ and 
consider the new generating function
$\tilde{f}=f+f'$.

\begin{proposition}
\label{eucau}
For a generic and small $f'$, 
$\tilde{f}$ has caustic diffeomorphic to a tricuspoid, the curve
shown in figure \ref{tricuspoid}.1 (see \cite{AGZV} for a definition of
tricuspoid).
\end{proposition}

\setlength{\unitlength}{1cm}
\begin{picture}(6,7)
\label{tricuspoid}
\thinlines
\put(0,3){\vector(1,0){6}}
\put(2,0){\vector(0,1){6}}
\put(5.5,3.2){$\displaystyle x_1$}
\put(2.2,5.7){$\displaystyle x_2$}
\qbezier(5,3)(2,3)(1,5)
\qbezier(5,3)(2,3)(1,1)
\qbezier(1,5)(2,3)(1,1)
\put(4.7,3.2){$\displaystyle A_3$}
\put(1,5.2){$\displaystyle A_3$}
\put(1,.6){$\displaystyle A_3$}
\put(1,3.2){$\displaystyle A_2$}
\put(3,3.5){$\displaystyle A_2$}
\put(3,2.2){$\displaystyle A_2$}
\end{picture}
\begin{center}
$Fig.~\ref{tricuspoid}.1:~The~tricuspoid$
\end{center}
Having only folds and cusps, the tricuspoid is stable. 
If $f'(y_1,y_2)=\frac{\epsilon}{2}
(ay_1^2+by_1y_2+cy_2^2)$ is a generic polynomial of degree $2$,
the critical locus turns out to be a circle in the $(y_1,y_2)$-plane
with centre 
$$C=\Big(-\frac{\epsilon}{4}(a-c),\frac{\epsilon}{4}b\Big)$$
and radius
$$\frac{\epsilon}{4}|a+c|$$
in this case the caustic is a tricuspoid and 
can be explicitly computed (see \cite{AGZV}).

\begin{proof}
By hypothesis, in the Whitney topology of $C^\infty(\R^2)$
$\tilde{f}$ lies in a small neighbourhood of $f$, so, if
$T$ is a tubular 
neighbourhood of $graph(f)\subset\R^2\times \R$, we can identify
$\tilde{f}$ with a section of $C^\infty(T)$ and find a deformation $\tilde{F}$
from $f$ to $\tilde{f}$. The Milnor number of $f$ is 4, thus a versal
deformation $F$ of $f$ has four parameters and can be written as
$F(y_1,y_2)=f(y_1,y_2)+a_0+a_1y_1+a_2y_2+a_3y_1^2$. 
By definition, any other deformation $G$ of $f$ is 
obtained from $F$ as $G(y,\lambda)=F(H(y,\lambda),
\Phi(\lambda))$, where $y=(y_1,y_2)$, $\lambda$ represents the parameters of
the deformation, $H$ is a family of diffeomorphisms parametrized by $\lambda$
and $\Phi$ is a smooth function of $\lambda$. Observe that $F$ is a
generating family of the elliptic umbilic in dimension 3 (in fact $f$
is the normal form of the singularities $D_4^-$):
it defines a generating function (see equation (\ref{eugenfun}))
$$\bar{f}(y_1,y_2,x_3)=y_1^3-y_1 y_2^2+x_3 y_1^2$$
and a Lagrangian map
\begin{equation*}
\left\{ \begin{aligned}
x_1 & = y_1^2-y_2^2+2x_3y_1\\
x_2 & = -2y_1y_2\\
y_3 & = -y_1^2
\end{aligned}
\right.
\end{equation*}
whose
caustic $K_F$, showed in figure \ref{pyramid}.2,
is the well known pyramid.\\\\
\setlength{\unitlength}{1cm}
\begin{picture}(8,5)
\label{pyramid}
\thinlines

\put(0,2){\vector(1,0){8}}
\put(4,0){\vector(0,1){5}}
\put(2,0){\vector(1,1){4}}
\put(7.7,2.1){$\displaystyle x_3$}
\put(4.1,4.7){$\displaystyle x_1$}
\put(5.7,4.2){$\displaystyle x_2$}

\qbezier(0.5,2.5)(1,2.5)(1.5,4)
\qbezier(1.5,4)(1.5,2)(2,1)
\qbezier(2,1)(1.5,2.5)(0.5,2.5)

\qbezier(4,2)(3,2)(1.5,4)
\qbezier(4,2)(3,2)(2,1)
\qbezier(4,2)(3,2.2)(0.5,2.5)

\qbezier(7.5,2.5)(7,2.5)(6.5,4)
\qbezier(6.5,4)(6.5,2)(6,1)
\qbezier(6,1)(6.5,2.5)(7.5,2.5)

\qbezier(4,2)(5,2)(6.5,4)
\qbezier(4,2)(5,2)(6,1)
\qbezier(4,2)(5,2.2)(7.5,2.5)

\put(0.6,3.3){$\displaystyle A_2$}
\put(5.7,2.6){$\displaystyle A_2$}
\put(2.2,3.3){$\displaystyle A_3$}
\put(7.6,2.5){$\displaystyle A_3$}
\put(6.5,1.6){$\displaystyle A_2$}
\put(4.8,1.3){$\displaystyle A_3$}
\put(2.1,1.8){$\displaystyle A_2$}
\put(3.8,2.2){$\displaystyle D_4^-$}
\end{picture}
\begin{center}
$Fig.~\ref{pyramid}.2:~The~pyramid$
\end{center}
Note that $f$ is recovered from $\bar{f}$ by setting $x_3=0$, so that
the caustic $K_f$ of $f$ can be 
identified with the intersection $K_F\cap\{x_3=0\}$
between the
pyramid and the plane 
$x_3=0$. 
Observe instead that the intersection $K_F\cap\{x_3=t\}$ is, for $t\neq0$,
a tricuspoid.
For $f'$ sufficiently small, $\tilde{F}$ is a
small deformation of the elliptic umbilic in dimension 3, and being
this stable, it follows that the caustic $K_{\tilde{F}}$
of $\tilde{F}$, in suitable coordinates $x_1'$, $x_2'$ and $x_3'$,
is still the pyramid. 
On the other hand, the versality of $F$
ensures the existence of a map $\Phi$ such that
$\Phi(x_3)=x_3'$ and relating, as explained, $\tilde{F}$ to $F$. 
Choosing $f'$ sufficiently small,
$\Phi$ will be enough
close to the identity, in the Whitney topology, to be
injective. Since 
$K_{\tilde{f}}=K_{{\tilde{F}}\cap\{x_3=0\}}$, it follows that the caustic
$K_{\tilde{f}}$ of $\tilde{f}$ 
is generically diffeomorphic to a tricuspoid.
\end{proof}
   
\subsection{The hyperbolic umbilic}
We refer to the generating function
\begin{equation} 
\label{hypumb}
f(y_1,y_2)=\frac{1}{3}(y_1^3+y_2^3)
\end{equation}
whose associated Lagrangian map is
\begin{equation*}
\left\{ \begin{aligned}
x_1 & =y_1^2\\
x_2 & =y_2^2
\end{aligned}
\right.
\end{equation*}
as the hyperbolic umbilic in dimension 2.
The critical locus is given by $y_1y_2=0$. 
The caustic is the set
$\{x_1x_2=0~:~x_1, x_2\geq0\}$.

\begin{proposition}
A generic 
small perturbation of the hyperbolic umbilic in dimension 2 has a caustic
diffeomorphic to the non-connected subset shown 
in figure \ref{hugra}.1.
\end{proposition}

\setlength{\unitlength}{1cm}
\begin{picture}(6,4)
\label{hugra}
\thinlines

\put(0,1.5){\vector(1,0){3.5}}
\put(1.5,0){\vector(0,1){3.5}}
\put(3.35,1.6){$\displaystyle x_1$}
\put(1.6,3.35){$\displaystyle x_2$}

\qbezier(2,2)(2.5,2.5)(2.5,3.5)
\qbezier(2,2)(2.5,2.5)(3.5,2.5)
\qbezier(.5,3.5)(.5,.5)(3.5,.5)

\put(1.9,1.7){$\displaystyle A_3$}
\put(2.9,2.1){$\displaystyle A_2$}
\put(1.9,2.8){$\displaystyle A_2$}
\put(.9,.9){$\displaystyle A_2$}
\end{picture}
\begin{center}
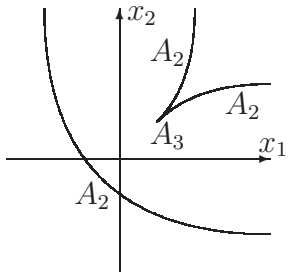
$Fig.~\ref{hugra}.1:~The~caustic~of~a~small~perturbation~of~the
~hyperbolic~umbilic$
\end{center}

\begin{proof}
The argument is the same as the one used in the proof of 
proposition \ref{eucau}.
\end{proof}

\subsection{Other singularities}
In dimension 2 we can consider other unstable germs of functions and try to
study how their caustics change, when slightly perturbed, by using their
versal deformations.
As seen in the previous subsections, 
being the elliptic and hyperbolic umbilics, and also the 
swallow-tail (see \cite{AGZV}), 
stable in dimension 3, the study of the generating functions 
(\ref{eugenfun2}) and (\ref{hypumb}), 
in dimension 2, was recovered from the analysis of the generating functions 
(\ref{eugenfun}) 
and (\ref{hugenfun}), in dimension 3,
by fixing one parameter. Instead, consider, for instance,  
the parabolic umbilic (see \cite{BL}), which is stable in dimension 4: 
it is necessary to fix two
parameters to recover the case of dimension 2 from the stable case
in dimension 4.
So, if
interested in some applications to mirror symmetry when the total space
of the fibration has complex dimension 2, 
it seems to be not so relevant to consider those unstable singularities,
such as the parabolic umbilic,
whose versal deformations define stable singularities in dimension
greater than 3: indeed, a generic orbit of Hamiltonian
equivalence containing a Lagrangian map exhibiting an unstable singularity,
such as the elliptic or hyperbolic umbilic, after a small perturbation, still
will contain a Lagrangian map with such singularity; this is no longer true
if the singularity is, for example, a parabolic umbilic.

\section{Gradient lines and their bifurcations}
\label{bifgradlin}
\subsection{Gradient lines of a Lagrangian map}
Given a Lagrangian map $L\hookrightarrow X\rightarrow B$ with generating 
function $f$, where we always assume $X=\R^4$ and $B=\R^2$, 
and fixed a metric on $X$, we define a family of
functions $f_x:\R^2\rightarrow\R$, parametrized by $x\in B$, as
$f_x(y)=f(y)-x\cdot y$, and consider a dynamical system on each
fibre $X_x=\R^2$ of $X$, over $x$, as follows:
\begin{equation}
\label{ds}
\frac{dy}{dt}=\nabla f_x
\end{equation}
where $\nabla$ is the gradient induced by the metric on $X$.

\begin{definition}
A curve $y:(a,b)\rightarrow X_x$, with $a$,
$b\in\R\cup\{+\infty,-\infty\}$, is called a gradient line if it is a
solution of (\ref{ds}).
\end{definition}

Note that the set of critical points of $\nabla f_x$ coincides with
the intersection $L\cap X_x$.

\begin{lemma}
If $x\notin K$, where $K$ is the caustic of $L$, then $f_x$ has only
non-degenerate critical points (in other words, $f_x$ is a Morse function).
\end{lemma}
\begin{proof}
If $y$ is a critical point of $f_x$, then $\nabla f_x(y)=\nabla
f(y)-x=0$. If $x\notin K$ then $Hf(y)=Hf_x(y)$ has maximal rank.
\end{proof}

Gradient vector fields share the following feature:

\begin{lemma}
\label{pt2}
If $f$ has only finitely many non-degenerate critical points, then $\nabla
f$ has finitely many fixed points all of which are hyperbolic and no
other periodic orbits.
\end{lemma}
\begin{proof}
See \cite{R}.
\end{proof}
 
The Morse index of a non-degenerate critical point $y$ of $f$ is defined as 
the number of negative eiegenvalues of the Hessian $Hf(y)$. In dimension 2, 
lemma \ref{pt2} implies that if $x$ does not belong to the
caustic, we expect as critical points of $\nabla f_x$ only unstable nodes, 
saddles and stable nodes, identified by Morse index respectively 
equal to 0, 1 and 2.

\subsection{Bifurcation points of a Lagrangian map}
\begin{definition}
A point $x\in B$ is a bifurcation point of $f$ 
if and only if $x\notin K$ and $\nabla f_x$ is not
Morse-Smale (see \cite{DM} or \cite{KH} or \cite{R}
for the definition of Morse-Smale vector field).
\end{definition}

\begin{corollary}
\label{corbif}
Let $x$ be a bifurcation point, then there exist two critical points $y_1$ and
$y_2$ of $\nabla f_x$ such that $W^u(y_1)$ and $W^s(y_2)$ do not intersect
transversely, where $W^u$ and $W^s$ denote respectively the unstable and 
stable manifold of critical points.
\end{corollary}
\begin{proof}
It is a direct consequence of the definition of Morse-Smale vector field and
of the fact that $x\notin K$ and that $\nabla f_x$ is a gradient vector field.
\end{proof}

\begin{remark}
\rm
Observe that for vector fields on 2-manifolds, the
Morse-Smale condition is equivalent to structural stability.
\end{remark}

In dimension 2, the critical points in corollary \ref{corbif} are saddles.
Since
the stable and unstable manifolds of a saddle 
are each the union of two of the four
separatrices of the saddle, a non-transversal intersection of
$W^u(y_1)$ and $W^s(y_2)$ means that
$y_1$ and $y_2$ have a common
separatrix, or, in other words, that there is a gradient line from
$y_1$ to $y_2$.
We call saddle-to-saddle separatrix such homoclinic orbit.

\begin{proposition}
A saddle-to-saddle separatrix is not structurally stable.
\end{proposition} 
\begin{proof}
See \cite{AN}
\end{proof}

Figure \ref{saddleconngra}.1 
shows the bifurcation given by a saddle-to-saddle
separatrix from $s_1$ to $s_2$.\\ 
\setlength{\unitlength}{1cm}
\begin{picture}(12,6)
\label{saddleconngra}
\thinlines

\qbezier(1,3)(1.5,3)(2,3)
\qbezier(2,3)(2.7,3)(2.7,5)
\qbezier(2,1)(2,3)(2,5)

\qbezier(2.3,1)(2.3,3)(3,3)
\qbezier(3,3)(4.5,3)(4,3)
\qbezier(3,1)(3,3)(3,5)

\qbezier(5,3)(6.5,3)(8,3)
\qbezier(7,1)(7,3)(7,5)
\qbezier(6,1)(6,3)(6,5)

\qbezier(9,3)(9.5,3)(10,3)
\qbezier(10,3)(10.7,3)(10.7,1)
\qbezier(10,1)(10,3)(10,5)

\qbezier(10.3,5)(10.3,3)(11,3)
\qbezier(11,3)(11.5,3)(12,3)
\qbezier(11,1)(11,3)(11,5)

\scriptsize
\put(1.7,3.1){$\displaystyle s_1$}
\put(3.1,3.1){$\displaystyle s_2$}
\put(5.7,3.1){$\displaystyle s_1$}
\put(7.1,3.1){$\displaystyle s_2$}
\put(9.7,3.1){$\displaystyle s_1$}
\put(11.1,3.1){$\displaystyle s_2$}
\tiny
\put(2.5,5){$W^u(s_1)$}
\put(2,.7){$W^s(s_2)$}
\put(6.1,3.1){$W^u(s_1)$}
\put(6.1,2.7){$W^s(s_2)$}
\put(10.5,.7){$W^u(s_1)$}
\put(10,5){$W^s(s_2)$}
\small
\put(1,1.5){$X_0$}
\put(5,1.5){$X_{bif}$}
\put(9,1.5){$X_1$}
\normalsize

\end{picture}
\begin{center}
$Fig.~\ref{saddleconngra}.1:~A~saddle-to-saddle~separatrix$
\end{center}
Observe that the structurally stable vector fields $X_0$ and $X_1$, 
though orbitally equivalent, 
are not orbitally equivalent under deformations, in the sense that, if
$\phi$ is the homeomorphism of the plane mapping the phase portrait of
$X_0$ to the phase portrait of $X_1$ and respecting the sense of the flow,
and if
$\Phi$ a homotopy between the identity and $\phi$, with parameter space 
$[0,1]$, then there exists $t\in(0,1)$ such that $\Phi(~,t)$ is not a
homeomorphism (in other words, as it is qualitatively evident, a
continuous deformation of the phase portrait of $X_0$ to the one of
$X_1$, respecting the direction of the flow, contains the phase
portrait of the unstable vector field $X_{bif}$).
Thus, for a generic family of vector fields exhibiting two saddles, near an
element having a saddle-to-saddle 
separatrix, there are two classes of vector fields
up to orbitally equivalence under deformations.

Denote by ${\cal M}(y_1,y_2)$ the moduli space of
unparametrized gradient lines from a critical point $y_1$ to a
critical point $y_2$.

\begin{proposition}
\label{morse}
If $\nabla f_x$ is Morse-Smale and ${\cal M}(y_1,y_2)\neq\emptyset$
then 
$$dim{\cal M}(y_1,y_2)=ind(y_1)-ind(y_2)-1$$
\end{proposition}
\begin{proof}
See for example \cite{DM}.
\end{proof}

This implies that 
gradient lines from $y_1$ to
$y_2$ exist generically only if the Morse index of $y_1$ is greater
than the Morse index of $y_2$, and they are stable.

Observe that, also for $x\in K$, the vector field $\nabla f_x$ is not
Morse-Smale: 
what happens is that the nature or the number of critical points of
$\nabla{f_x}$ change. These bifurcations are called local bifurcations,
because it is enough to study the
vector field in a neighbourhood of the degenerate bifurcating critical
points. Instead, those bifurcations involving a lack of transversality
between the stable and unstable manifolds of two critical points, as
in the case of a saddle-to-saddle separatrix, are called
global, since involving global properties of the flow of the 
field $\nabla f_x$.

\begin{definition}
The bifurcation locus ${\cal B}$ of $f$ is the set of bifurcation points of
$f$. The diagram containing the caustic $K$ and 
the bifurcation locus ${\cal B}$ of $f$ in $B=\R^2$ is the
bifurcation diagram of $f$.
\end{definition}

Each point $x$ of a bifurcation diagram gives information about critical
points and existence of saddle-to-saddle separatrices of $\nabla
f_x$. Far from the caustic $K$, the vector field $\nabla f_x$ exhibits
a certain number of saddles $s_1(x)$,...,$s_n(x)$, so 
we can define components ${\cal B}_{i,j}$ of the bifurcation locus
${\cal B}$ 
as the 
set of points $x$ such that $\nabla f_x$ exhibits a gradient line 
$\gamma_{s_i(x)s_j(x)}$ from $s_i(x)$ to $s_j(x)$. 

\begin{proposition}
\label{codimbif}
Far from $K$ and from other components of ${\cal B}$, 
${\cal B}_{i,j}$, if non-empty, is an immersed
submanifold of codimension 1.
\end{proposition}
\begin{proof}
Let $S(s_i(x_0))$ and $S(s_j(x_0))$ be respectively the separatrices of 
$s_i(x_0)$ and
of $s_j(x_0)$ which intersect, at $x_0\in{\cal B}_{i,j}$, in the gradient
line $\gamma_{s_i(x_0)s_j(x_0)}$. 
Consider
a neighbourhood $N(x_0)$ of $x_0$ such that $N(x_0)$
does not intersect $K$ or other components of ${\cal B}$
different from ${\cal B}_{i,j}$, then, for all $x\in N(x_0)$,
the vector field $\nabla f_x$, if structurally stable, 
belongs to two distinct classes ${\cal V}_1$ and 
${\cal V}_2$ up to orbital equivalence under deformations.
Define $\psi:N(x_0)\rightarrow\R$ as  
\begin{displaymath}
\psi(x)=\left\{ \begin{array}{ll}
dist(S(s_i(x)),S(s_j(x)))^2 & x\in{\cal V}_1\\
-dist(S(s_i(x)),S(s_j(x)))^2 & x\in{\cal V}_2
\end{array} \right.
\end{displaymath}
Note that $\psi$ is smooth everywhere, because the family $f_x$ depends 
smoothly on $x$, and that,
if non-empty, ${\cal B}_{i,j}=\psi^{-1}(0)$ 
(this is true because $N(x_0)$
does not intersect $K$ or other components of ${\cal B}$
different from ${\cal B}_{i,j}$: in fact if a saddle $s_k(x)$, with 
$k\neq i,j$, were a limit point of both $S(s_i(x))$ and $S(s_j(x))$,
then $\psi(x)=0$ 
though there is no gradient line from $s_i(x)$ to
$s_j(x)$). Generically, $\psi$ is a Morse function, 
thus ${\cal B}$ is an immersed submanifold of $N(x_0)$; 
moreover,
far from its critical points, $\psi$ is transversal to
$0\in\R$, so ${\cal B}_{i,j}$ is a 
submanifold of $N(x_0)$ of codimension 1. 
\end{proof}

In a similar way we can define subsets ${\cal B}_{(i,j),(k,l)}$ of
${\cal B}$ as the set of points $x\in B$ where $\nabla f_x$ exhibits both the
exceptional gradient lines $\gamma_{s_i(x)s_j(x)}$ and
$\gamma_{s_k(x)s_l(x)}$.

\begin{corollary}
Far from $K$ and from other components
of ${\cal B}\setminus({\cal B}_{i,j}\cup{\cal B}_{k,l})$, 
${\cal B}_{(i,j),(k,l)}$, if non-empty,
is an immersed submanifold of codimension 2.  
\end{corollary}
\begin{proof}
Define, in a neighbourhood $N(x_0)$ of a point 
$x\in{\cal B}_{(i,j),(k,l)}$ which is far from 
$K$ and other components of 
${\cal B}\setminus({\cal B}_{i,j}\cup{\cal B}_{k,l})$, a function
$\psi:B\rightarrow \R^2$ as $\psi(x)=(\psi_{ij},\psi_{kl})$, where
$\psi_{ij}$ and $\psi_{kl}$ are as in the proof of proposition
\ref{codimbif}, and note that ${\cal B}_{(i,j),(k,l)}=\psi^{-1}(0)$ and
$\{0\}$ has codimension 2 in $\R^2$.
\end{proof}

For a generic $f$, ${\cal B}_{(i,j),(k,l)}={\cal B}_{i,j}\cap 
{\cal B}_{k,l}$. 
It is clear that three exceptional gradient lines in
the same phase portrait is
a bifurcation of codimension greater than 2, 
so, generically, it does not occur in
dimension 2. Therefore, ${\cal B}$ can be decomposed into 
strata ${\cal B}_{i,j}$ and ${\cal B}_{(i,j),(k,l)}$, whose codimension
is respectively 1 and 2. 

Whether exceptional gradient lines appear or not, or in other
words, whether the subsets ${\cal B}_{i,j}$ and ${\cal B}_{(i,j),(k,l)}$ are
non-empty, it depends on the family of vector fields. 
We will analyze those cases which we need to study the bifurcation diagram
of the cusp and of the elliptic umbilic. 
Since we are dealing with the family $\nabla
f_x$, we may assume that all the elements of the family
are gradient vector fields, in order to avoid troubles with periodic orbits.

\subsection{A family of vector fields with two saddles}
\label{saddlesaddle}
Consider a 2-parameters family of gradient vector fields $X_x$ exhibiting two
saddles $s_1$ and $s_2$. Consider a
structurally stable element $X_s$ of the family. By definition of stability,
there exists a neighbourhood $U$ of $s$ such that, for every $t\in U$,
$X_t$ is conjugated to $X_s$. On the boundary $\partial U$ 
of $U$ we can expect to
meet a bifurcation point $b$, where $X_b$ presents a
saddle-to-saddle separatrix. 

\begin{proposition}
A point $b\in{\cal B}\cap\partial U$ can belong to either 
${\cal B}_{1,2}$ or ${\cal B}_{2,1}$. 
\end{proposition}
\begin{proof}
The saddle-to saddle separatrix of $X_b$ can be given
by either $W^u(s_1)\cap W^s(s_2)$ or by $W^s(s_1)\cap W^u(s_2)$, 
which means that
$b$ belongs respectively to ${\cal B}_{1,2}$ or to 
${\cal B}_{2,1}$.
\end{proof}

By proposition \ref{codimbif}, both ${\cal B}_{1,2}$ and ${\cal B}_{2,1}$ have
codimension 1.

\begin{lemma}
\label{saddlesaddleinters}
${\cal B}_{1,2}\cap{\cal B}_{2,1}=\emptyset$. The intersection of 
two components 
${\cal B}_{ij}^1$ and ${\cal B}_{ij}^2$ of
${\cal B}_{ij}$ is non-empty provided the saddle-to-saddle separatrices
$\gamma_{ij}^1$ and $\gamma_{ij}^2$, appearing respectively at points
of ${\cal B}_{ij}^1$ and ${\cal B}_{ij}^2$, are obtained as
intersection of the same pair of 
separatrices of $s_1$ and $s_2$.
\end{lemma}
\begin{proof}
Let $b\in{\cal B}_{1,2}\cap{\cal B}_{2,1}$,
then $X_b$ exhibits two saddles and, between them, two
saddle-to-saddle separatrices with opposite directions;
consider in $\R^2$ a close curve $C$ containing $s_1$ and $s_2$, then the
Poincar\'e index $ind_{\cal P}(C)$ of $C$ would be equal to -1, while, on the
other hand, $ind_{\cal P}(C)=ind_{\cal P}(s_1)+ind_{\cal P}(s_2)=-2$.
This proves the first statement.

For $t\in{\cal B}_{ij}^1\cap{\cal B}_{ij}^2$, suppose the gradient lines
$\gamma_{ij}^1$ and $\gamma_{ij}^2$ of $X_t$
are obtained as intersection of different pairs of separatrices of the saddles,
then the phase portrait of $X_t$
exhibits two exceptional gradient lines between the two saddles, 
giving a contraddiction as shown in the first part of the proof.
Otherwise, no contraddiction arises at $t$, since only one
saddle-to-saddle separatrix appears in the phase portrait of $X_t$. Moreover, 
if $\alpha$, $\beta$,
$\gamma$ and $\delta$ denote the four subsets determined in $\R^2$ by
${\cal B}_{ij}^1$ and ${\cal B}_{ij}^2$,\\
\setlength{\unitlength}{1cm}
\begin{picture}(4,4)
\label{try}
\thinlines
\qbezier(1,1)(2,2)(3,3)
\qbezier(1,3)(2,2)(3,1)

\put(.4,.7){${\cal B}_{ij}^1$}
\put(.4,3.1){${\cal B}_{ij}^1$}
\put(1.3,1.9){$\alpha$}
\put(2.4,1.9){$\gamma$}
\put(1.9,2.4){$\beta$}
\put(1.9,1.4){$\delta$}
\end{picture}
\begin{center}
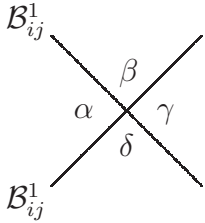
$Fig.~\ref{try}.1:~The~intersection~of~two~components~of~{\cal B}_{ij}$
\end{center}
then the two classes of orbitally equivalent vector fields 
under deformation are given
by $x\in\alpha\cup\gamma$ and $x\in\beta\cup\delta$ (see figure \ref{try}.1) .
\end{proof}

\subsection{A family of vector fields with two saddles and one node}
\label{twosaddlesnode}
Consider a 2-parameters family of gradient vector fields $X_x$ with two
saddles $s_1$ and $s_2$ and an unstable node $n$. 
We want to understand which kind of bifurcations, that is, which kind of
saddle-to-saddle separatrices, the family can exhibit.
Consider, if existing, a
structurally stable element $X_s$ of the family, such that its phase
portrait contains both the gradient lines $\gamma_{ns_1}$ and
$\gamma_{ns_2}$, from $n$ to respectively $s_1$ and $s_2$. 
Structural stability ensures the existence of an open connected
neighbourhood $U_1$ of $s$ in $\R^2$ such that, for every $t\in U_1$, the phase
portrait of $X_t$ is orbitally equivalent to the phase portrait of $X_s$.
Among such neighbourhoods of $s$ we can assume $U_1$ to be maximal.

\begin{proposition}
\label{snprop}
If ${\cal B}\cap\partial U_1\neq\emptyset$,
a point $t\in{\cal B}\cap\partial U_1$ belongs to either 
${\cal B}_{1,2}$ or ${\cal B}_{2,1}$: the saddle-to-saddle separatrix
$\gamma_{s_is_j}$ is obtained as intersection of $\gamma_{ns_j}$ with
one of the two components of $W^u(s_i)$; at $t$, where $\gamma_{s_is_j}$
appears, $\gamma_{ns_j}$ breaks. 
\end{proposition}
\begin{proof}
For $t\in{\cal B}\cap\partial U_1$,
the saddle-to-saddle separatrices which can be exhibited in the phase portrait
of $X_t$ 
are $\gamma_{s_1s_2}$ and $\gamma_{s_2s_1}$,
implying that $t$ belongs respectively to
${\cal B}_{1,2}$ and ${\cal B}_{2,1}$. Consider, for
example, $\gamma_{s_2s_1}$ ($\gamma_{s_1s_2}$ can be treated similarly): 
$\gamma_{s_2s_1}=W^u(s_2)\cap W^s(s_1)$; as shown  
in figure \ref{2snode}.1\\
\setlength{\unitlength}{1cm}
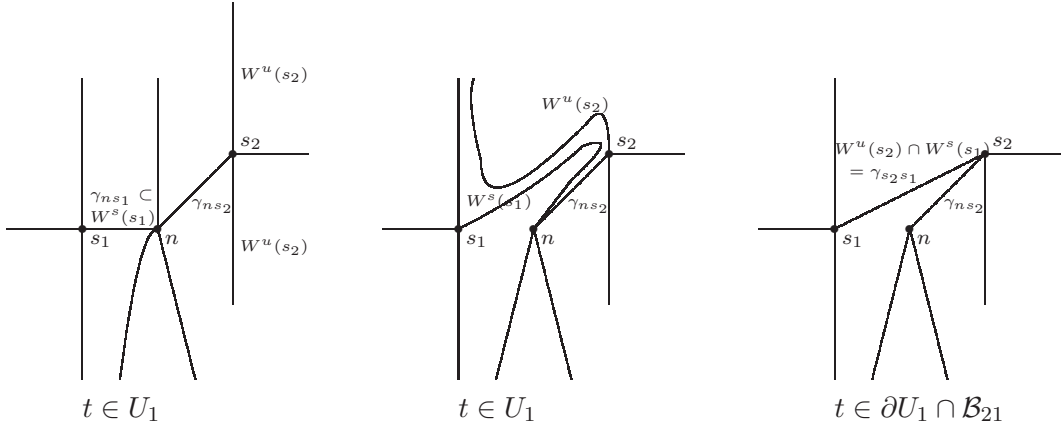
\begin{picture}(12,6.5)
\label{2snode}
\thinlines

\put(1,3){\circle*{.1}}
\put(2,3){\circle*{.1}}
\put(3,4){\circle*{.1}}

\put(6,3){\circle*{.1}}
\put(7,3){\circle*{.1}}
\put(8,4){\circle*{.1}}

\put(11,3){\circle*{.1}}
\put(12,3){\circle*{.1}}
\put(13,4){\circle*{.1}}

\scriptsize
\put(1.1,2.8){$\displaystyle s_1$}
\put(2.1,2.8){$\displaystyle n$}
\put(3.1,4.1){$\displaystyle s_2$}

\put(6.1,2.8){$\displaystyle s_1$}
\put(7.1,2.8){$\displaystyle n$}
\put(8.1,4.1){$\displaystyle s_2$}

\put(11.1,2.8){$\displaystyle s_1$}
\put(12.1,2.8){$\displaystyle n$}
\put(13.1,4.1){$\displaystyle s_2$}

\tiny
\put(2.45,3.3){$\gamma_{ns_2}$}
\put(1.1,3.4){$\gamma_{ns_1}\subset$}
\put(1.1,3.1){$W^s(s_1)$}
\put(3.1,2.7){$W^u(s_2)$}
\put(3.1,5){$W^u(s_2)$}

\put(7.45,3.3){$\gamma_{ns_2}$}
\put(12.45,3.3){$\gamma_{ns_2}$}

\put(11.2,3.7){$=\gamma_{s_2s_1}$}
\put(11,4){$W^u(s_2)\cap W^s(s_1)$}

\put(6.1,3.3){$W^s(s_1)$}
\put(7.1,4.6){$W^u(s_2)$}

\small
\put(1,.5){$t\in U_1$}
\put(6,.5){$t\in U_1$}
\put(11,.5){$t\in \partial U_1\cap{\cal B}_{21}$}

\normalsize

\qbezier(1,3)(1.5,3)(2,3)
\qbezier(1,1)(1,3)(1,5)
\qbezier(1,3)(0.5,3)(0,3)

\qbezier(3,4)(2.5,3.5)(2,3)
\qbezier(3,2)(3,4)(3,6)
\qbezier(3,4)(3.5,4)(4,4)

\qbezier(2,3)(2,4)(2,5)
\qbezier(2,3)(1.75,3)(1.5,1)
\qbezier(2,3)(2.25,2)(2.5,1)

\qbezier(6,3)(7,3.5)(7.65,4.1)
\qbezier(7,3)(7.5,3.65)(7.5,3.65)
\qbezier(7.65,4.1)(8.2,4.3)(7.5,3.65)
\qbezier(6,1)(6,3)(6,5)
\qbezier(6,3)(5.5,3)(5,3)

\qbezier(8,4)(7.5,3.5)(7,3)
\qbezier(8,2)(8,3)(8,4)
\qbezier(8,4)(8,4.8)(7.7,4.4)
\qbezier(7.7,4.4)(6.3,3)(6.3,3.9)
\qbezier(6.3,3.9)(6.1,4.7)(6.2,5)
\qbezier(8,4)(8.5,4)(9,4)

\qbezier(7,3)(6.75,2)(6.5,1)
\qbezier(7,3)(7.25,2)(7.5,1)

\qbezier(11,3)(12,3.5)(13,4)
\qbezier(11,1)(11,3)(11,5)
\qbezier(11,3)(10.5,3)(10,3)

\qbezier(13,4)(12.5,3.5)(12,3)
\qbezier(13,2)(13,3)(13,4)
\qbezier(13,4)(13.5,4)(14,4)

\qbezier(12,3)(11.75,2)(11.5,1)
\qbezier(12,3)(12.25,2)(12.5,1)

\end{picture}\\
\setlength{\unitlength}{1cm}
\begin{picture}(12,5)
\thinlines

\put(1,3){\circle*{.1}}
\put(2,3){\circle*{.1}}
\put(3,4){\circle*{.1}}

\scriptsize
\put(1.1,2.8){$\displaystyle s_1$}
\put(2.1,2.8){$\displaystyle n$}
\put(3.1,4.1){$\displaystyle s_2$}

\tiny
\put(.1,1.5){$W^u(s_1)$}
\put(2.1,5){$W^s(s_1)$}
\put(3.4,3.7){$W^s(s_2)$}
\put(1.85,3.9){$W^u(s_2)$}
\put(2.45,3.3){$\gamma_{ns_2}$}

\small
\put(1,.5){$t\in U_2$}
\normalsize

\qbezier(1,3)(2,3.5)(3,5)
\qbezier(1,1)(1,3)(1,5)
\qbezier(1,3)(0.5,3)(0,3)

\qbezier(3,4)(2.5,3.5)(2,3)
\qbezier(3,2)(3,3)(3,4)
\qbezier(3,4)(1.8,3.7)(1.2,1.5)
\qbezier(3,4)(3.5,4)(4,4)

\qbezier(2,3)(1.75,3)(1.5,1)
\qbezier(2,3)(2.25,2)(2.5,1)

\end{picture}
\begin{center}
$Fig.~\ref{2snode}.1:~The~bifurcation~from~U_1~to~U_2$
\end{center}
the only connected component of $W^s(s_1)$ which can intersect $W^u(s_2)$ is 
$\gamma_{ns_1}$, while 
both the unstable separatrices of $s_2$, 
the connected component of $W^u(s_2)$,
can intersect $W^s(s_1)$. 
This implies that
when, at the bifurcation point $t$, $\gamma_{s_2s_1}$ appears, 
$\gamma_{ns_1}$ breaks.
\end{proof}

Let $U_2$ be a (maximal)
open connected subset such that $X_t$ is structurally stable for all $t\in U_2$
and ${\cal B}_{ij}\cap\partial U_1\cap\partial U_2\neq\emptyset$.

\begin{proposition}
\label{snprop2}
For all $t\in U_2$,
the phase portrait of $X_t$ does not exhibit the gradient
line $\gamma_{ns_j}$.
At a point $t\in\partial U_2\cap{\cal B}$, 
two pairs of separatrices, one of $s_i$ and one of $s_j$, can intersect in
a saddle-to-saddle separatrix:
in one case, described in
proposition \ref{snprop} and shown in figure \ref{2snode}.1, 
$t$ belongs to 
${\cal B}_{ij}$, moreover, after the bifurcation, the line $\gamma_{ns_j}$
appears in the phase portrait of $X_t$; 
in the other case, analyzed 
in subsection
\ref{saddlesaddle} for a family of vector fields with two saddles,
$t$ belongs to ${\cal B}_{ji}$, moreover after the bifurcation no gradient 
line $\gamma_{ns_j}$ appears in the phase portrait of $X_t$ (see figure
\ref{2snode2}.2).  
\end{proposition}
\begin{proof}
That for $t\in U_2$ the phase portrait of $X_t$ does not contain
$\gamma_{ns_j}$ is a consequence of proposition \ref{snprop}. Setting
for simplicity $i=2$ and $j=1$ as in the proof of proposition
\ref{snprop},
the two pair of separatrices of $s_1$ and $s_2$ that can intersect
at $t\in\partial U_2\cap{\cal B}$ are shown in figure \ref{2snode}.1 
and in figure \ref{2snode2}.2. In the first case,
the saddle-to-saddle separatrix is $\gamma_{s_2s_1}=W^u(s_2)\cap W^s(s_1)$,
so $t\in{\cal B}_{21}$, and after the bifurcation $W^s(s_1)=\gamma_{ns_1}$.
In the second case, shown in figure \ref{2snode2}.2, 
$\gamma_{s_1s_2}=W^u(s_1)\cap W^s(s_2)$,
so $t\in{\cal B}_{12}$. Observe moreover that two  
separatrices of $s_2$, among those 
not intersecting with $W^u(s_1)$ in $\gamma_{s_1s_2}$,
determines in $\R^2$ two disjoint subsets, one containing $n$ and one 
containing $s_1$, which implies that after the bifurcation $\gamma_{ns_1}$
does not appear in the phase portrait of $X_t$
(so the bifurcation is of the type described
in subsection 
\ref{saddlesaddle}).\\
\setlength{\unitlength}{1cm}
\begin{picture}(12,5.5)
\label{2snode2}
\thinlines

\put(1,3){\circle*{.1}}
\put(2,3){\circle*{.1}}
\put(3,4){\circle*{.1}}

\scriptsize
\put(1.1,2.8){$\displaystyle s_1$}
\put(2.1,2.8){$\displaystyle n$}
\put(3.1,4.1){$\displaystyle s_2$}

\tiny
\put(.1,1.5){$W^u(s_1)$}
\put(3.4,3.7){$W^s(s_2)$}
\put(2.45,3.3){$\gamma_{ns_2}$}

\small
\put(1,.5){$t\in U_2$}
\normalsize

\qbezier(1,3)(2,3.5)(3,5)
\qbezier(1,1)(1,3)(1,5)
\qbezier(1,3)(0.5,3)(0,3)

\qbezier(3,4)(2.5,3.5)(2,3)
\qbezier(3,2)(3,3)(3,4)
\qbezier(3,4)(1.8,3.7)(1.2,1.5)
\qbezier(3,4)(3.5,4)(4,4)

\qbezier(2,3)(1.75,3)(1.5,1)
\qbezier(2,3)(2.25,2)(2.5,1)

\thinlines

\put(6,3){\circle*{.1}}
\put(7,3){\circle*{.1}}
\put(8,4){\circle*{.1}}

\scriptsize
\put(6.1,2.8){$\displaystyle s_1$}
\put(7.1,2.8){$\displaystyle n$}
\put(8.1,3.9){$\displaystyle s_2$}

\tiny
\put(.1,1.5){$W^u(s_1)$}
\put(3.4,3.7){$W^s(s_2)$}
\put(7.45,3.3){$\gamma_{ns_2}$}

\put(6.65,4){$=\gamma_{s_1s_2}$}
\put(6.7,4.2){$W^u(s_1)\cap W^s(s_2)$}

\small
\put(6,.5){$t\in\partial U_2\cap{\cal B}_{12}$}
\normalsize

\qbezier(6,3)(6.5,3.5)(6.5,5)
\qbezier(6,3)(6,4)(6,5)
\qbezier(6,3)(5.5,3)(5,3)
\qbezier(6,3)(6.5,3)(7,3.5)

\qbezier(8,4)(7.5,3.5)(7,3)
\qbezier(8,2)(8,3)(8,4)
\qbezier(8,4)(6.8,3.7)(6.2,1.5)
\qbezier(8,4)(7.5,4)(7,3.5)

\qbezier(7,3)(6.75,3)(6.5,1)
\qbezier(7,3)(7.25,2)(7.5,1)

\thinlines

\put(11,3){\circle*{.1}}
\put(12,3){\circle*{.1}}
\put(13,4){\circle*{.1}}

\scriptsize
\put(11.1,2.8){$\displaystyle s_1$}
\put(12.1,2.8){$\displaystyle n$}
\put(13.1,3.8){$\displaystyle s_2$}

\tiny
\put(12.1,4.8){$W^u(s_1)$}
\put(10,2.8){$W^s(s_1)$}
\put(12.1,4.1){$W^s(s_2)$}
\put(13.1,2.2){$W^u(s_2)$}
\put(12.45,3.3){$\gamma_{ns_2}$}


\small
\put(11,.5){$t\in U_3$}
\normalsize

\qbezier(11,3)(11.5,3.5)(11.5,5)
\qbezier(11,3)(11,4)(11,5)
\qbezier(11,3)(10.5,3)(10,3)
\qbezier(11,3)(11.5,3)(12,5)

\qbezier(13,4)(12.5,3.5)(12,3)
\qbezier(13,2)(13,3)(13,4)
\qbezier(13,4)(11.8,3.7)(11.2,1.5)
\qbezier(13,4)(11.7,4)(11,1.5)

\qbezier(12,3)(11.75,3)(11.5,1)
\qbezier(12,3)(12.25,2)(12.5,1)

\end{picture}
\begin{center}
$Fig.~\ref{2snode2}.2:~The~bifurcation~from~U_2~to~U_3$
\end{center}
\end{proof}

Let $U_3$ be a (maximal)
open connected subset such that $X_t$ is structurally stable for all $t\in U_3$
and ${\cal B}_{ji}\cap\partial U_2\cap\partial U_3\neq\emptyset$.

\begin{proposition}
\label{snprop3}
For all $t\in U_3$,
the phase portrait of $X_t$ does not exhibit the gradient
line $\gamma_{ns_j}$.
At a point $t\in\partial U_3\cap{\cal B}$, 
two pairs of separatrices, one of $s_i$ and one of $s_j$, can intersect in
a saddle-to-saddle separatrix:
in one case, described in
proposition \ref{snprop2} and shown in figure \ref{2snode}.2, 
$t$ belongs to 
${\cal B}_{ji}$, and after the bifurcation 
the phase portrait of $X_t$ does not exhibit
the line $\gamma_{ns_j}$;
in the other, shown in figure \ref{2snode3}.3, $t$ belongs to 
${\cal B}_{ij}$, and after the bifurcation also 
the phase portrait of $X_t$ does not exhibit
the line $\gamma_{ns_j}$;
both bifurcations are of the type
analyzed 
in subsection
\ref{saddlesaddle} for a family of vector fields with two saddles.
\end{proposition}
\begin{proof}
The proposition can be proved as done for proposition \ref{snprop2}.
The saddle-to-saddle separatrices at $t\in\partial U_3\cap{\cal B}$, in the
two cases,
are shown respectively 
in figure \ref{2snode2}.2 and in figure \ref{2snode3}.3 below.\\
\setlength{\unitlength}{1cm}
\begin{picture}(12,5.5)
\label{2snode3}
\thinlines

\put(1,3){\circle*{.1}}
\put(2,3){\circle*{.1}}
\put(3,4){\circle*{.1}}

\scriptsize
\put(1.1,2.8){$\displaystyle s_1$}
\put(2.1,2.8){$\displaystyle n$}
\put(3.1,3.8){$\displaystyle s_2$}

\tiny
\put(0,2.8){$W^s(s_1)$}
\put(3.1,2.2){$W^u(s_2)$}
\put(2.45,3.3){$\gamma_{ns_2}$}


\small
\put(1,.5){$t\in U_3$}
\normalsize

\qbezier(1,3)(1.5,3.5)(1.5,5)
\qbezier(1,3)(1,4)(1,5)
\qbezier(1,3)(0.5,3)(0,3)
\qbezier(1,3)(1.5,3)(2,5)

\qbezier(3,4)(2.5,3.5)(2,3)
\qbezier(3,2)(3,3)(3,4)
\qbezier(3,4)(1.8,3.7)(1.2,1.5)
\qbezier(3,4)(1.7,4)(1,1.5)

\qbezier(2,3)(1.75,3)(1.5,1)
\qbezier(2,3)(2.25,2)(2.5,1)

\thinlines

\put(6,3){\circle*{.1}}
\put(7,3){\circle*{.1}}
\put(8,4){\circle*{.1}}

\scriptsize
\put(5.8,2.8){$\displaystyle s_1$}
\put(7.1,2.8){$\displaystyle n$}
\put(8.1,3.8){$\displaystyle s_2$}

\tiny
\put(7.45,3.3){$\gamma_{ns_2}$}

\put(6.65,4.1){$=\gamma_{s_1s_2}$}
\put(6.7,4.3){$W^u(s_1)\cap W^s(s_2)$}

\small
\put(6,.5){$t\in\partial U_3\cap{\cal B}_{21}$}
\normalsize

\qbezier(6,3)(6.2,3.5)(6.2,5)
\qbezier(6,3)(6,4)(6,5)
\qbezier(6,3)(6.5,3.4)(7,3.7)
\qbezier(6,3)(6.3,3)(6.5,5)

\qbezier(8,4)(7.5,3.5)(7,3)
\qbezier(7,3.7)(7.7,4.1)(8,4)
\qbezier(8,4)(6.8,3.7)(6.2,1.5)
\qbezier(8,4)(6.7,4)(6,1.5)

\qbezier(7,3)(6.75,3)(6.5,1)
\qbezier(7,3)(7.25,2)(7.5,1)

\thinlines

\put(11,3){\circle*{.1}}
\put(12,3){\circle*{.1}}
\put(13,4){\circle*{.1}}

\scriptsize
\put(10.8,2.8){$\displaystyle s_1$}
\put(12.1,2.8){$\displaystyle n$}
\put(13.1,3.8){$\displaystyle s_2$}

\tiny
\put(12.7,3.55){$\gamma_{ns_2}$}
\put(12,4.35){$\gamma_{ns_1}$}


\small
\put(11,.5){$t\in\ U_4$}
\normalsize

\qbezier(11,3)(11.2,3.5)(11.2,5)
\qbezier(11,3)(11,4)(11,5)
\qbezier(11,3)(11.3,3.1)(11.8,4)
\qbezier(13.3,4.1)(11.9,4.4)(11.8,4)
\qbezier(13.3,4.1)(14.5,3.3)(12,3)
\qbezier(11,3)(11.3,3)(11.5,5)

\qbezier(13,4)(11.6,4.1)(10.8,1.5)
\qbezier(12,3)(12.5,3.5)(13,4)
\qbezier(13,4)(11.8,3.7)(11.2,1.5)
\qbezier(13,4)(11.7,4)(11,1.5)

\qbezier(12,3)(11.75,3)(11.5,1)
\qbezier(12,3)(12.25,2)(12.5,1)

\end{picture}
\begin{center}
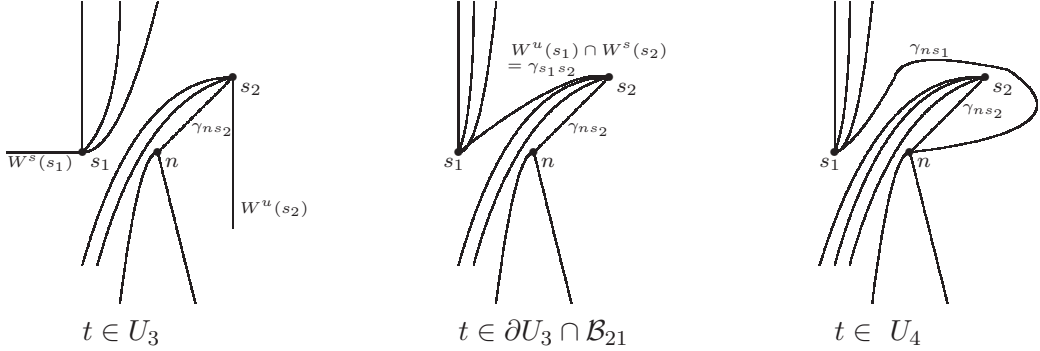
$Fig.~\ref{2snode3}.3:~The~bifurcation~from~U_3~to~U_4$
\end{center}
\end{proof}

Let $U_4$ be a (maximal)
open connected subset such that $X_t$ is structurally stable for all $t\in U_4$
and ${\cal B}_{ji}\cap\partial U_3\cap\partial U_4\neq\emptyset$.

\begin{proposition}
\label{snprop4}
For $t\in U_4$, the
phase portrait of $X_t$ does exhibit the gradient
line $\gamma_{ns_j}$ (see figure \ref{2snode3}.3).
\end{proposition}
\begin{proof}
It follows from proposition \ref{snprop3}. 
\end{proof}

Note that $\gamma_{ns_j}$ has a different winding around $s_i$ for 
$t\in U_1$ and $t\in U_4$.
 
\begin{lemma}
\label{ssnintersec}
${\cal B}_{1,2}\cap {\cal B}_{2,1}=\emptyset$. 
The intersection of 
two components 
${\cal B}_{ij}^1$ and ${\cal B}_{ij}^2$ of
${\cal B}_{ij}$ is non-empty provided the saddle-to-saddle separatrices
$\gamma_{ij}^1$ and $\gamma_{ij}^2$ 
of ${\cal B}_{ij}^1$ and ${\cal B}_{ij}^2$ are obtained as
intersection of the same pair of 
separatrices of $s_1$ and $s_2$.
\end{lemma}
\begin{proof}
See the proof
of lemma
\ref{saddlesaddleinters}.
\end{proof}

\subsection{A family of vector fields with three saddles and one node}
\label{3sn}
Consider a 2-parameter family of gradient vector fields $X_x$ exhibiting three
saddles $s_1$, $s_2$ and $s_3$ and an unstable node $n$. 
Consider a
structurally stable element $X_s$ of the family, such that its phase
portrait contains all the gradient lines $\gamma_{ns_i}$ for $i=1,2,3$, 
then there exists a
neighbourhood $U_1$ of $s$ such that, for every $t\in U_1$, the phase
portrait of $X_t$ is orbitally equivalent to the one of $X_s$ and so it
shows the same qualitative features. We can assume $U_1$ to be maximal
among such neighbourhoods of $s$.
At a point $t\in\partial U_1$, the vector field $X_t$
can exhibit one among the saddle-to-saddle separatrices 
$\gamma_{s_is_j}$,
implying $t\in{\cal B}_{ij}$. 
We have $\gamma_{s_is_j}=W^u(s_i)\cap W^s(s_j)$: in this case
$W^s(s_j)=\gamma_{ns_j}$, and so at $t$, where $\gamma_{s_is_j}$ appears, 
$\gamma_{ns_j}$
breaks. Moreover, unlike what described for a family of vector
fields exhibiting two saddles and a node, the choice of the component of
$W^u(s_i)$ is fixed by the presence of $\gamma_{ns_k}$, for $k\neq i,j$.
In figure \ref{2snode4}.1 the case of the saddle-to-saddle separatrix
$\gamma_{s_2s_1}$ is outlined.\\
\setlength{\unitlength}{1cm}
\begin{picture}(14,4.5)
\label{2snode4}

\thinlines
\put(1,2){\circle*{.1}}
\put(2.7,2.7){\circle*{.1}}
\put(2.7,1.3){\circle*{.1}}
\put(2,2){\circle*{.1}}

\scriptsize
\put(1.1,1.8){$\displaystyle s_1$}
\put(2.2,2){$\displaystyle n$}
\put(2.9,2.7){$\displaystyle s_2$}
\put(2.9,1.3){$\displaystyle s_3$}
\normalsize

\qbezier(2,2)(1.5,2)(1,2)
\qbezier(2,2)(2.5,2.5)(2.7,2.7)
\qbezier(2,2)(2.5,1.5)(2.7,1.3)

\qbezier(1,2)(1,3)(1,3.5)
\qbezier(1,2)(1,1)(1,.5)
\qbezier(1,2)(.5,2)(0,2)

\qbezier(2.7,2.7)(3,3)(3.5,3.5)
\qbezier(2.7,2.7)(2.4,3)(2,3.4)
\qbezier(2.7,2.7)(3,2.4)(3.5,2.2)

\qbezier(2.7,1.3)(3,1)(3.5,.6)
\qbezier(2.7,1.3)(3,1.6)(3.5,1.8)
\qbezier(2.7,1.3)(2.4,1)(2,.6)

\put(1,.1){$t\in U_1$}
\thinlines
\put(5,2){\circle*{.1}}
\put(6.7,2.7){\circle*{.1}}
\put(6.7,1.3){\circle*{.1}}
\put(6,2){\circle*{.1}}

\scriptsize
\put(5.1,1.8){$\displaystyle s_1$}
\put(6.2,2){$\displaystyle n$}
\put(6.9,2.7){$\displaystyle s_2$}
\put(6.9,1.3){$\displaystyle s_3$}
\normalsize

\qbezier(6,2)(6.5,2.5)(6.7,2.7)
\qbezier(6,2)(6.5,1.5)(6.7,1.3)

\qbezier(5,2)(5,3)(5,3.5)
\qbezier(5,2)(5,1)(5,.5)
\qbezier(5,2)(6,2.41)(6.7,2.7)
\qbezier(5,2)(4.5,2)(4,2)

\qbezier(6.7,2.7)(7,3)(7.5,3.5)
\qbezier(6.7,2.7)(7,2.4)(7.5,2.2)

\qbezier(6.7,1.3)(7,1)(7.5,.6)
\qbezier(6.7,1.3)(7,1.6)(7.5,1.8)
\qbezier(6.7,1.3)(6.4,1)(6,.6)

\put(5,.1){$t\in\partial U_1\cap{\cal B}_{ij}$}
\thinlines
\put(9,2){\circle*{.1}}
\put(10.7,2.7){\circle*{.1}}
\put(10.7,1.3){\circle*{.1}}
\put(10,2){\circle*{.1}}

\scriptsize
\put(9.1,1.8){$\displaystyle s_1$}
\put(10.2,2){$\displaystyle n$}
\put(10.9,2.7){$\displaystyle s_2$}
\put(10.9,1.3){$\displaystyle s_3$}
\normalsize

\qbezier(10,3.2)(9.5,2.5)(9,2)
\qbezier(10,2)(10.5,2.5)(10.7,2.7)
\qbezier(10,2)(10.5,1.5)(10.7,1.3)

\qbezier(9,2)(9,3)(9,3.5)
\qbezier(9,2)(9,1)(9,.5)
\qbezier(9,2)(8.5,2)(8,2)

\qbezier(10.7,2.7)(11,3)(11.5,3.5)
\qbezier(10.7,2.7)(9.3,2.2)(9.3,.7)
\qbezier(10.7,2.7)(11,2.4)(11.5,2.2)

\qbezier(10.7,1.3)(11,1)(11.5,.6)
\qbezier(10.7,1.3)(11,1.6)(11.5,1.8)
\qbezier(10.7,1.3)(10.4,1)(10,.6)

\put(9,.1){$t\in U_2$}


\end{picture}
\begin{center}
$Fig.~\ref{2snode4}.1:~The~bifurcation~from~U_1~to~U_2$
\end{center}
Let $U_2$ be a (maximal)
open connected subset such that $X_t$ is structurally stable for all $t\in U_2$
and ${\cal B}_{ij}\cap\partial U_1\cap\partial U_2\neq\emptyset$.
For $t\in U_2$, the
phase portrait of $X_t$ contains the gradient line $\gamma_{ns_i}$ and
$\gamma_{ns_k}$ for $k\neq i,j$,  but not 
$\gamma_{ns_j}$. Note that two of the separatrices of $s_i$ determine 
two
subsets
in the plane, one containing the node $n$ and the saddle $s_k$, $k\neq i,j$,
and the other containing the saddle
$s_j$.
This implies that $W^u(s_j)\cap W^s(s_k)=W^s(s_j)\cap
W^u(s_k)=\emptyset$, thus $\partial U_2\cap{\cal B}_{jk}=\partial U_2\cap
{\cal B}_{kj}=\emptyset$. 
Instead, $\partial U_2$
can intersect ${\cal B}_{i,j}$, as just described, and also 
${\cal B}_{j,i}$, ${\cal B}_{i,k}$ or ${\cal B}_{k,i}$.
As to the intersection with 
${\cal B}_{j,i}$, it holds what already outlined and shown in figures
of subsection
\ref{twosaddlesnode}: indeed, a separatrix of $s_j$ divides the plane
into two subsets, one containing the saddle $s_k$ and its separatrices and
one containing the pair of separatrices of $s_i$ and $s_j$ intersecting
in the saddle-to-saddle separatrix $\gamma_{s_j,s_i}$; moreover, when
$\gamma_{s_j,s_i}$ breaks, the phase portrait of $X_t$ does not exhibit
$\gamma_{n,s_i}$.
Instead, for what concerns ${\cal B}_{i,k}$,
observe that the gradient line
$\gamma_{ns_j}$ does not appear in the phase portrait of $X_t$, 
so it follows that
$W^s(s_k)=\gamma_{ns_k}$ can intersect both the separatrices defining 
$W^u(s_i)$, as shown in figures \ref{3bifurcations}.2 and 
\ref{3bifurcations2}.3. 
Analogous considerations hold for ${\cal B}_{k,i}$.\\
\setlength{\unitlength}{1cm}
\begin{picture}(14,4.5)
\label{3bifurcations}

\thinlines
\put(1,2){\circle*{.1}}
\put(2.7,2.7){\circle*{.1}}
\put(2.7,1.3){\circle*{.1}}
\put(2,2){\circle*{.1}}

\scriptsize
\put(1.1,1.8){$\displaystyle s_1$}
\put(2.2,2){$\displaystyle n$}
\put(2.9,2.7){$\displaystyle s_2$}
\put(2.9,1.3){$\displaystyle s_3$}
\normalsize

\qbezier(2,2)(2.5,2.5)(2.7,2.7)
\qbezier(2,2)(2.5,1.5)(2.7,1.3)

\qbezier(1,2)(1,3)(1,3.5)
\qbezier(1,2)(1,1)(1,.5)
\qbezier(1,2)(.5,2)(0,2)
\qbezier(2,3.2)(1.5,2.5)(1,2)

\qbezier(2.7,2.7)(3,3)(3.5,3.5)
\qbezier(2.7,2.7)(3,2.4)(3.5,2.2)
\qbezier(2.7,2.7)(1.3,2.2)(1.3,.7)

\qbezier(2.7,1.3)(3,1)(3.5,.6)
\qbezier(2.7,1.3)(3,1.6)(3.5,1.8)
\qbezier(2.7,1.3)(2.4,1)(2,.6)

\put(1,.1){$t\in U_2$}

\thinlines
\put(5,2){\circle*{.1}}
\put(6.7,2.7){\circle*{.1}}
\put(6.7,1.3){\circle*{.1}}
\put(6,2){\circle*{.1}}

\scriptsize
\put(5.1,1.8){$\displaystyle s_1$}
\put(6.2,2){$\displaystyle n$}
\put(6.9,2.7){$\displaystyle s_2$}
\put(6.9,1.3){$\displaystyle s_3$}
\normalsize

\qbezier(6,2)(6.5,2.5)(6.7,2.7)

\qbezier(5,2)(5,3)(5,3.5)
\qbezier(5,2)(5,1)(5,.5)
\qbezier(5,2)(4.5,2)(4,2)
\qbezier(6,3.2)(5.5,2.5)(5,2)

\qbezier(6.7,2.7)(7,3)(7.5,3.5)
\qbezier(6.7,2.7)(5.3,2.2)(5.3,.7)

\qbezier(6.7,1.3)(7,1)(7.5,.6)
\qbezier(6.7,1.3)(7,1.6)(7.5,1.8)
\qbezier(6.7,1.3)(6.4,1)(6,.6)
\qbezier(6.7,2.7)(6.7,2)(6.7,1.3)

\put(5,.1){$t\in\partial U_2\cap{\cal B}_{ik}$}

\thinlines
\put(9,2){\circle*{.1}}
\put(10.7,2.7){\circle*{.1}}
\put(10.7,1.3){\circle*{.1}}
\put(10,2){\circle*{.1}}

\scriptsize
\put(9.1,1.8){$\displaystyle s_1$}
\put(10.2,2){$\displaystyle n$}
\put(10.9,2.7){$\displaystyle s_2$}
\put(10.9,1.3){$\displaystyle s_3$}
\normalsize

\qbezier(10,3.2)(9.5,2.5)(9,2)
\qbezier(10,2)(10.5,2.5)(10.7,2.7)

\qbezier(9,2)(9,3)(9,3.5)
\qbezier(9,2)(9,1)(9,.5)
\qbezier(9,2)(8.5,2)(8,2)

\qbezier(10.7,2.7)(11,3)(11.5,3.5)
\qbezier(10.7,2.7)(9.3,2.2)(9.3,.7)
\qbezier(10.7,2.7)(10.4,1.2)(9.7,.8)

\qbezier(10.7,1.3)(11,1)(11.5,.6)
\qbezier(10.7,1.3)(11,1.6)(11.5,1.8)
\qbezier(10.7,1.3)(10.4,1)(10,.6)
\qbezier(10.7,1.3)(11,1.9)(11.5,2.3)

\put(9,.1){$t\in U_3$}


\end{picture}
\begin{center}
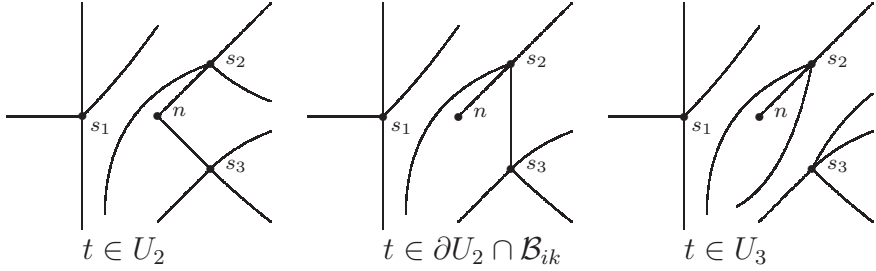
$Fig.~\ref{3bifurcations}.2:~The~bifurcation~from~U_2~to~U_3~(1^{st}~case)$
\end{center}
\setlength{\unitlength}{1cm}
\begin{picture}(14,4.5)
\label{3bifurcations2}

\thinlines
\put(1,2){\circle*{.1}}
\put(2.7,2.7){\circle*{.1}}
\put(2.7,1.3){\circle*{.1}}
\put(2,2){\circle*{.1}}

\scriptsize
\put(1.1,1.8){$\displaystyle s_1$}
\put(2.2,2){$\displaystyle n$}
\put(2.9,2.7){$\displaystyle s_2$}
\put(2.9,1.3){$\displaystyle s_3$}
\normalsize

\qbezier(2,2)(2.5,2.5)(2.7,2.7)
\qbezier(2,2)(2.5,1.5)(2.7,1.3)

\qbezier(1,2)(1,3)(1,3.5)
\qbezier(1,2)(1,1)(1,.5)
\qbezier(1,2)(.5,2)(0,2)
\qbezier(2,3.2)(1.5,2.5)(1,2)

\qbezier(2.7,2.7)(3,3)(3.5,3.5)
\qbezier(2.7,2.7)(3,2.4)(3.5,2.2)
\qbezier(2.7,2.7)(1.3,2.2)(1.3,.7)

\qbezier(2.7,1.3)(3,1)(3.5,.6)
\qbezier(2.7,1.3)(3,1.6)(3.5,1.8)
\qbezier(2.7,1.3)(2.4,1)(2,.6)

\put(1,.1){$t\in U_2$}

\thinlines
\put(5,2){\circle*{.1}}
\put(6.7,2.7){\circle*{.1}}
\put(6.7,1.3){\circle*{.1}}
\put(6,2){\circle*{.1}}

\scriptsize
\put(5.1,1.8){$\displaystyle s_1$}
\put(6.2,2){$\displaystyle n$}
\put(6.9,2.7){$\displaystyle s_2$}
\put(6.9,1.3){$\displaystyle s_3$}
\normalsize

\qbezier(6,2)(6.5,2.5)(6.7,2.7)

\qbezier(5,2)(5,3)(5,3.5)
\qbezier(5,2)(5,1)(5,.5)
\qbezier(5,2)(4.5,2)(4,2)
\qbezier(6,3.2)(5.5,2.5)(5,2)

\qbezier(6.7,2.7)(7,3)(7.5,3.5)
\qbezier(6.7,2.7)(7,2.4)(7.5,2.2)
\qbezier(6.7,2.7)(4.5,2)(6.7,1.3)

\qbezier(6.7,1.3)(7,1)(7.5,.6)
\qbezier(6.7,1.3)(7,1.6)(7.5,1.8)
\qbezier(6.7,1.3)(6.4,1)(6,.6)

\put(5,.1){$t\in\partial U_2\cap{\cal B}_{ki}$}

\thinlines
\put(9,2){\circle*{.1}}
\put(10.7,2.7){\circle*{.1}}
\put(10.7,1.3){\circle*{.1}}
\put(10,2){\circle*{.1}}

\scriptsize
\put(9.1,1.8){$\displaystyle s_1$}
\put(10.2,2){$\displaystyle n$}
\put(10.9,2.7){$\displaystyle s_2$}
\put(10.9,1.3){$\displaystyle s_3$}
\normalsize

\qbezier(10,3.2)(9.5,2.5)(9,2)
\qbezier(10,2)(10.5,2.5)(10.7,2.7)

\qbezier(9,2)(9,3)(9,3.5)
\qbezier(9,2)(9,1)(9,.5)
\qbezier(9,2)(8.5,2)(8,2)

\qbezier(10.7,2.7)(11,3)(11.5,3.5)
\qbezier(10.7,2.7)(8.5,1.5)(11.5,2)
\qbezier(10.7,2.7)(11,2.4)(11.5,2.2)

\qbezier(10.7,1.3)(11,1)(11.5,.6)
\qbezier(10.7,1.3)(11,1.6)(11.5,1.8)
\qbezier(10.7,1.3)(10.4,1)(10,.6)
\qbezier(10.3,3)(8.5,2)(10.7,1.3)

\put(9,.1){$t\in U_3$}

\put(12,0){$\displaystyle Fig.~\ref{3bifurcations2}.3$}

\end{picture}
\begin{center}
$Fig.~\ref{3bifurcations}.3:~The~bifurcation~from~U_2~to~U_3~(2^{nd}~case)$
\end{center}
From figure \ref{3bifurcations}.2 
we see that for $t\in U_3$ 
no intersection is
possible between the stable and unstable manifolds of $s_1$ and $s_3$,
so $\partial U_3$ can intersect only ${\cal B}_{12}$,  
${\cal B}_{21}$,  ${\cal B}_{23}$ or ${\cal B}_{32}$.
From figure \ref{3bifurcations2}.3 we see instead that for
$t\in U_3$
only  ${\cal B}_{13}$, ${\cal B}_{31}$, ${\cal B}_{23}$ and 
${\cal B}_{32}$ can intersect $\partial U_3$.

We can resume what said in the following, very general, proposition:

\begin{proposition}
\label{allowed}
The gradient line $\gamma_{s_is_j}$ can
appear if a
component of $W^u(s_i)$ can intersect a component of $W^s(s_j)$.
Whether this is possible and which components can actually intersect
depends on the separatrices of the third saddle $s_k$ and on the
gradient lines $\gamma_{ns_l}$, $l=1,2,3$, appearing in the phase portrait.
\end{proposition}

As to intersection of bifurcation lines, we will analyze some cases
when studying the bifurcation locus of perturbations of the elliptic
umbilic. The caustic, as we will see, imposes further constraints on
the possible intersections.

\subsection{A family of vector fields with saddles and nodes}
\label{ssgencase}
In general,
to a family $X_x$, $x\in\R^2$, of planar vector fields, having nodes 
$n_1$, ..., $n_\alpha$, and saddles $s_1$, ...,$s_\beta$, 
we can associate a bifurcation diagram given
by the bifurcation locus ${\cal B}$ with its components ${\cal B}_{ij}$
and ${\cal B}_{(ij),(k,l)}$. 
A proposition similar to \ref{allowed} can be formulated also in this case, 
explaining where, depending on the phase portrait of $X_x$,  
the gradient lines $\gamma_{s_is_j}$ 
appear.

\begin{proposition}
\label{gencase}
Let $U\subset\R^2$ be 
an open connected subset such that for every $x\in U$, the vector fields
$X_x$ are structurally stable and orbitally equivalent. 
Then $\partial U$ can intersect ${\cal B}_{ij}$ if and only if $W^u(s_i)$
and $W^s(s_j)$ lie in the same connected component determined in the
phase portrait of $X_x$ by:\\
- the separatrices of the remaining saddles of $X_x$\\
- the
gradient lines $\gamma_{n_ks_j}$ or $\gamma_{s_jn_k}$, 
respectively if $n_k$ is an unstable or stable node, for 
$k=1,...,\alpha$ and $j=1,...,\beta$. 
\end{proposition}

As to intersection of bifurcation lines, as said, generically,
${\cal B}_{(i,j),(j,k)}={\cal B}_{(i,j)}\cap{\cal B}_{(j,k)}$. 
A necessary condition for ${\cal B}_{(i,j),(k,l)}$ to be non-empty is that 
there exists a vector field whose phase portrait can exhibit both the 
exceptional gradient lines $\gamma_{s_is_j}$ and $\gamma_{s_ks_l}$. Moreover,
in a neighbourhood $N(t)$ of $t\in{\cal B}_{(i,j),(k,l)}$, 
the saddle-to-saddle separatrices
generically break, and the phase portrait of $\nabla f_x$, for $x\in N(t)$, 
must be recovered continuosly
from the phase portrait of $\nabla f_t$.
In particular, if $t\in{\cal B}_{(i,j),(j,k)}$,
the exceptional gradient lines $\gamma_{s_is_j}$ and $\gamma_{s_js_k}$,
can also break in $N(t)$ in such a way to form
the saddle-to-saddle separatrix $\gamma_{s_is_k}$, implying
$t\in\overline{{\cal B}_{(i,k)}}$.
As to intersection of
bifurcation lines,
we will consider some examples
when studying perturbations of the elliptic umbilic.

\begin{definition}
Given a bifurcation diagram, 
expressions as ``the diagram is allowed'' or ``permitted'' will be used to
mean that there exists a continuous family of 
planar vector fields providing the given bifurcation diagram.
\end{definition}

For example, a bifurcation diagram, such that
${\cal B}_{ij}\cap{\cal B}_{ji}\neq\emptyset$, is not allowed.
 
Observe also that a bifurcation diagram contains information about
the existence of non-generic gradient lines but no information about
the number and nature of critical points. 

\subsection{Families of gradient fields}
Given an allowed diagram, the problem is now to understand, at least in
those example we are concerned with, when
there exists a family of vector fields of the form $\nabla f_x$,
where $f_x(y)=f(y)-x\cdot y$ with $f:\R^2\rightarrow\R$.
The first
step is to construct a family of gradient vector fields whose
bifurcation diagram is the given one, and then to look for a family with
the required dependence from the parameter.

%
\begin{lemma}
\label{integration1}
Suppose 
that in an open simply connected subset $U\subset\R^2$ the phase portrait of a
vector field $X$ does not exhibit
any critical points or periodic orbits, then there exists a function
$f$ on $U$ such that $\nabla f$ is orbitally equivalent to $X$.
\end{lemma}
\begin{proof}
Consider the distribution of vector fields $\{X_x\}_{x\in U}$. Since
$X(x)\neq0$ for $x\in U$, we can choose an orthogonal
distribution $\{X^\perp_x\}_{x\in U}$, and since $U$ is simply
connected we can suppose this distribution to be smooth. The
hypothesis of Frobenius theorem are satisfied, so through every point
$x\in U$ it passes a unique curve integrating the distribution
$\{X^\perp_x\}_{x\in U}$. In a
neighbourhood $V$ of any point $x\in U$ a function $f_V$ can be defined
having the curves integrating $\{X^\perp_x\}_{x\in U}$ in $V$ as level
curves. Since the integral curves of a gradient vector field cross the
level sets of its potential ortogonally at points which are not fixed 
points, it follows that
the phase portrait in $V$ of $\nabla
f_V$ coincides with the phase portrait of $X_{|V}$. Since $U$ is
simply connected a function $f$ having the property required can be
defined on the whole $U$.
\end{proof}

 Unfortunately, given a phase portrait exhibiting only saddles and
nodes, it does not necessarily exist a potential $f$ whose gradient field
$\nabla f$ exhibits that phase portrait. We can state the following lemma,
but not its converse:

\begin{lemma}
\label{crpgrf}
If $p$ is a local maximum, minimum or saddle of $f$, then
$p$ is respectively an unstable node, a stable node or a saddle of 
$\nabla f$.
\end{lemma}
\begin{proof}
It is a consequence of the definition of node and saddle of a vector field.
\end{proof}

\begin{lemma}
\label{fundamental}
Suppose $f:\R^2\rightarrow \R$ is a function such that $\nabla
f$ is a Morse-Smale vector field exhibiting two saddles, then
there exists a function\\ 
$F:\R^2\times[-1,1]\rightarrow\R$ such that
$F(~,-1)=f$, $\nabla F(~,t)$ is a Morse-Smale vector field
exhibiting two saddles for every $t\neq0$, and $\nabla F(~,0)$ has a
saddle-to-saddle separatrix.
\end{lemma}
\begin{proof}
Lemma \ref{crpgrf} tells which behaviour the level curves of $F(~,t)$
have in a neighbourhood of saddles, and moreover we know
these level curves are orthogonal to the separatrices
of the saddles. In figure \ref{potential1}.1, the saddles $s_1^t$ and $s_2^t$ 
of $\nabla F(~,t)$,
their separatrices, denoted by $a_i^t$ and $b_i^t$, $i=1,2$,
and some of the relevant level curves of $F(~,t)$ (the red lines), 
for $-1\leq t<0$,
are shown.\\
\setlength{\unitlength}{1cm}
\begin{picture}(7,4.5)
\label{potential1}
\thinlines

\put(2,2){\circle*{.1}}
\put(4,2){\circle*{.1}}

\qbezier(0,2)(1,2)(2,2)
\qbezier(2,0)(2,2)(2,4)

\qbezier(6,2)(5,2)(4,2)
\qbezier(4,0)(4,2)(4,4)

\qbezier(2,2)(3.5,2)(3.5,4)
\qbezier(4,2)(2.5,2)(2.5,0)

\thicklines
\color{red}

\qbezier(2,2)(1,1)(0.5,0.5)
\qbezier(2,2)(1,3)(0.5,3.5)
\qbezier(2,2)(2.5,2.5)(2.8,4)
\qbezier(2,2)(2.5,1.5)(3,0)

\qbezier(4,2)(5,1)(5.5,0.5)
\qbezier(4,2)(5,3)(5.5,3.5)
\qbezier(4,2)(3.5,1.5)(3.2,0)
\qbezier(4,2)(3.5,2.5)(3,4)

\qbezier(2.3,2.2)(2.35,2.05)(2.4,1.95)
\qbezier(2.5,2.3)(2.7,2.1)(2.8,2)
\qbezier(2.7,2.4)(2.9,2.2)(3.1,2.1)

\qbezier(3.7,1.8)(3.65,1.95)(3.6,2.05)
\qbezier(3.5,1.7)(3.3,1.9)(3.2,2)
\qbezier(3.3,1.6)(3.1,1.8)(2.9,1.9)
\normalcolor

\scriptsize
\put(1.4,2.1){$s_1^t$}
\put(4.3,2.1){$s_2^t$}
\put(2.8,2.5){$a_1^t$}
\put(3,1.4){$a_2^t$}
\put(2.1,1){$b_1$}
\put(4.1,3){$b_1$}

\color{red}
\normalcolor
\normalsize

\end{picture}
\begin{center}
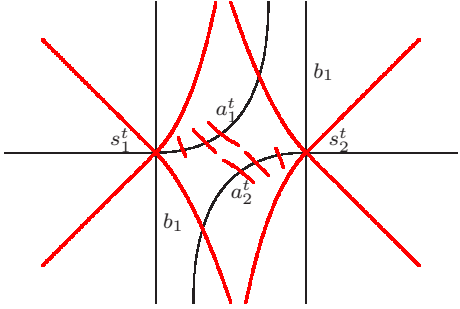
$Fig.~\ref{potential1}.1:~Phase~portrait~of~\nabla F(~,t)~and~level~curves~of~
F(~,t)$
$for~-1\leq t<0$
\end{center}
In figure \ref{potential2}.2, the phase portrait of $\nabla F(~,t)$ and some of
the relevant level curves of $F(~,t)$ for respectively $t=0$ and
$0<t\leq1$, are shown. The saddle-to saddle separatrix is denoted by $a^0$.\\
\setlength{\unitlength}{1cm}
\begin{picture}(13,4.5)
\label{potential2}
\thinlines

\put(2,2){\circle*{.1}}
\put(4,2){\circle*{.1}}

\qbezier(0,2)(1,2)(2,2)
\qbezier(2,0)(2,2)(2,4)

\qbezier(6,2)(5,2)(4,2)
\qbezier(4,0)(4,2)(4,4)

\qbezier(2,2)(3,2)(4,2)

\thicklines
\color{red}

\qbezier(2,2)(1,1)(0.5,0.5)
\qbezier(2,2)(1,3)(0.5,3.5)
\qbezier(2,2)(2.5,2.5)(2.8,4)
\qbezier(2,2)(2.5,1.5)(2.8,0)

\qbezier(4,2)(5,1)(5.5,0.5)
\qbezier(4,2)(5,3)(5.5,3.5)
\qbezier(4,2)(3.5,1.5)(3.2,0)
\qbezier(4,2)(3.5,2.5)(3.2,4)

\qbezier(2.7,2.3)(2.7,2)(2.7,1.7)
\qbezier(3.3,2.3)(3.3,2)(3.3,1.7)
\qbezier(2.3,2.15)(2.3,2)(2.3,1.85)
\qbezier(3.7,2.15)(3.7,2)(3.7,1.85)

\normalcolor

\scriptsize
\put(1.4,2.1){$s_1^0$}
\put(4.3,2.1){$s_2^0$}
\put(2.9,2.1){$a^0$}
\put(2.1,1){$b_1^0$}
\put(4.1,3){$b_2^0$}
\normalsize
\thinlines

\put(9,2){\circle*{.1}}
\put(11,2){\circle*{.1}}

\qbezier(7,2)(8,2)(9,2)
\qbezier(9,0)(9,2)(9,4)

\qbezier(13,2)(12,2)(11,2)
\qbezier(11,0)(11,2)(11,4)

\qbezier(9,2)(10.5,2)(10.5,0)
\qbezier(11,2)(9.5,2)(9.5,4)

\thicklines
\color{red}

\qbezier(9,2)(8,1)(7.5,0.5)
\qbezier(9,2)(8,3)(7.5,3.5)
\qbezier(9,2)(9.5,2.5)(10,4)
\qbezier(9,2)(9.5,1.5)(9.8,0)

\qbezier(11,2)(12,1)(12.5,0.5)
\qbezier(11,2)(12,3)(12.5,3.5)
\qbezier(11,2)(10.5,1.5)(10,0)
\qbezier(11,2)(10.5,2.5)(10.2,4)

\qbezier(10.7,2.2)(10.65,2.05)(10.6,1.95)
\qbezier(10.5,2.3)(10.3,2.1)(10.2,2)
\qbezier(10.3,2.4)(10.1,2.2)(9.9,2.1)

\qbezier(9.3,1.8)(9.35,1.95)(9.4,2.05)
\qbezier(9.5,1.7)(9.7,1.9)(9.8,2)
\qbezier(9.7,1.6)(9.9,1.8)(10.1,1.9)
\normalcolor

\scriptsize
\put(8.4,2.1){$s_1^t$}
\put(11.3,2.1){$s_2^t$}
\put(9.9,2.5){$a_1^t$}
\put(10.2,1.5){$a_2^t$}
\put(9.1,1){$b_1^t$}
\put(11.1,3){$b_2^t$}
\normalsize
\end{picture}
\begin{center}
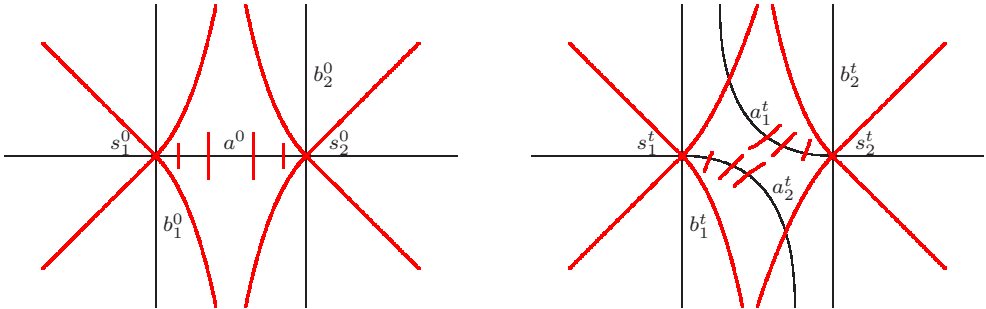
$Fig.~\ref{potential1}.2:~Phase~portrait~of~\nabla F(~,t)~and~level~curves~of~
F(~,t)$
$for~t=0~and~0<t\leq1~respectively$
\end{center}
To construct the functions $F(~,t)$,
we first choose two points in $\R^2$, in whose 
neighbourhoods we define, 
according to lemma \ref{crpgrf}, $F(~,t)$
in such a way
that these points are saddles; then,
we set $F(x,t)=f(x)$ for every $x\in\R^2\setminus A$ and
$t\in[-1,1]$, where $A\subset\R^2$ is a
neighbourhood, shown in figure \ref{potential3}.3, 
of the line chosen as the saddle-to-saddle separatrix
of $F(0,t)$.\\
\setlength{\unitlength}{1cm}
\begin{picture}(13,4.5)
\label{potential3}
\thinlines

\put(2,2){\circle*{.1}}
\put(4,2){\circle*{.1}}

\qbezier(0,2)(1,2)(2,2)
\qbezier(2,0)(2,2)(2,4)

\qbezier(6,2)(5,2)(4,2)
\qbezier(4,0)(4,2)(4,4)

\qbezier(2,2)(3,2)(4,2)

\thicklines
\color{red}

\qbezier(2,2)(2.5,2.5)(2.8,4)
\qbezier(2,2)(2.5,1.5)(2.8,0)

\qbezier(4,2)(3.5,1.5)(3.2,0)
\qbezier(4,2)(3.5,2.5)(3.2,4)

\qbezier(2.8,0)(3,0)(3.2,0)
\qbezier(2.8,4)(3,4)(3.2,4)

\normalcolor

\scriptsize
\put(1.6,2){$s_1$}
\put(4.2,2){$s_2$}
\put(2.9,2.1){$a$}
\put(2.1,1){$b_1$}
\put(4.1,3){$b_2$}
\put(2.85,1){$A$}
\normalsize

\end{picture}
\begin{center}
$Fig.~\ref{potential3}.3:~The~subset~A$
\end{center}
We define the level curves of $F(~,t)$ in a
neighbourhood of $a_i^t$ as the fibres of the normal bundle to $a_i^t$ and
extend then $F(~,t)$ to the whole $A$ (see figure \ref{2potential1}.4): 
this means to require
that the derivatives of $F(~,t)$ in the direction normal to the
separatrices $a_i^t$ vanish at any point of the separatrices $a_i^t$
$$\frac{\partial F(~,t)}{\partial (a_i^t)^\perp}(a_i^t(s))=0$$\\
\setlength{\unitlength}{1cm}
\begin{picture}(7,5.5)
\label{2potential1}
\thinlines

\put(2,2){\circle*{.1}}
\put(4,2){\circle*{.1}}

\qbezier(2,2)(3.5,2)(3.5,4)
\qbezier(4,2)(2.5,2)(2.5,0)

\thicklines
\color{red}

\qbezier(2,2)(2.5,2.5)(2.8,4)
\qbezier(2,2)(2.5,1.5)(2.8,0)

\qbezier(4,2)(3.5,1.5)(3.2,0)
\qbezier(4,2)(3.5,2.5)(3.2,4)

\qbezier(2.8,0)(3,0)(3.2,0)
\qbezier(2.8,4)(3,4)(3.2,4)

\qbezier(2.3,2.2)(2.35,2.05)(2.4,1.95)
\qbezier(2.5,2.3)(2.7,2.1)(2.8,2)
\qbezier(2.7,2.4)(2.9,2.2)(3.1,2.1)

\qbezier(3.7,1.8)(3.65,1.95)(3.6,2.05)
\qbezier(3.5,1.7)(3.3,1.9)(3.2,2)
\qbezier(3.3,1.6)(3.1,1.8)(2.9,1.9)
\normalcolor

\scriptsize
\put(1.6,2){$s_1^t$}
\put(4.2,2){$s_2^t$}
\put(2.8,2.5){$a_1^t$}
\put(3,1.4){$a_2^t$}
\put(2.85,1){$A$}

\color{red}
\normalcolor
\normalsize

\thinlines

\put(2,2){\circle*{.1}}
\put(4,2){\circle*{.1}}

\qbezier(6,2)(7,2)(8,2)

\thicklines
\color{red}

\qbezier(6,2)(6.5,2.5)(6.8,4)
\qbezier(6,2)(6.5,1.5)(6.8,0)

\qbezier(8,2)(7.5,1.5)(7.2,0)
\qbezier(8,2)(7.5,2.5)(7.2,4)

\qbezier(6.8,0)(7,0)(7.2,0)
\qbezier(6.8,4)(7,4)(7.2,4)

\qbezier(6.7,2.3)(6.7,2)(6.7,1.7)
\qbezier(7.3,2.3)(7.3,2)(7.3,1.7)
\qbezier(6.3,2.15)(6.3,2)(6.3,1.85)
\qbezier(7.7,2.15)(7.7,2)(7.7,1.85)

\normalcolor

\scriptsize
\put(5.6,2){$s_1^0$}
\put(8.2,2){$s_2^0$}
\put(6.9,2.1){$a^0$}
\put(6.85,1){$A$}
\normalsize

\normalsize
\thinlines

\put(10,2){\circle*{.1}}
\put(12,2){\circle*{.1}}

\qbezier(10,2)(11.5,2)(11.5,0)
\qbezier(12,2)(10.5,2)(10.5,4)

\thicklines
\color{red}

\qbezier(11.7,2.2)(11.65,2.05)(11.6,1.95)
\qbezier(11.5,2.3)(11.3,2.1)(11.2,2)
\qbezier(11.3,2.4)(11.1,2.2)(10.9,2.1)

\qbezier(10.3,1.8)(10.35,1.95)(10.4,2.05)
\qbezier(10.5,1.7)(10.7,1.9)(10.8,2)
\qbezier(10.7,1.6)(10.9,1.8)(11.1,1.9)

\qbezier(10,2)(10.5,2.5)(10.8,4)
\qbezier(10,2)(10.5,1.5)(10.8,0)

\qbezier(12,2)(11.5,1.5)(11.2,0)
\qbezier(12,2)(11.5,2.5)(11.2,4)

\qbezier(10.8,0)(11,0)(11.2,0)
\qbezier(10.8,4)(11,4)(11.2,4)

\normalcolor

\scriptsize
\put(9.6,2){$s_1^t$}
\put(12.2,2){$s_2^t$}
\put(10.9,2.5){$a_1^t$}
\put(11.2,1.5){$a_2^t$}
\put(10.85,1){$A$}
\normalsize

\put(2,4.5){$-1\leq t<0$}
\put(6.5,4.5){$t=0$}
\put(10,4.5){$0<t\leq1$}
\end{picture}
\begin{center}
$Fig.~\ref{2potential1}.4:~The~construction~of~F$
\end{center}
By construction, $\nabla F(~,t)$ has the required properties.
\end{proof} 

Observe that, for each $t$, the conditions defining the function $F(~,t)$
concerns only those points which we choose as critical points or 
belonging to separatrices 
in $A$.

\begin{corollary}
\label{fundamentalcor}
Suppose $f:\R^2\rightarrow \R$ is a function such that $\nabla
f$ is a Morse-Smale vector field exhibiting two saddles, then
there exists a function\\ 
$F:\R^2\times D^2\rightarrow\R$, where
$D^2=\{t_1^2+t_2^2\leq1\}\subset\R^2$, such that
$F(~,0,-1)=f$, $\nabla F(~,t_1,t_2)$ is a Morse-Smale vector field
exhibiting two saddles for every $t_2\neq0$, and $\nabla F(~,t_1,0)$ has a
saddle-to-saddle separatrix.
\end{corollary}
\begin{proof}
The requirement that the family $F(~,t_1,t_2)$ exhibits a
saddle-to-saddle separatrix along the subset $\{t_2=0\}\subset D^2$ is
compatible with what said about the dimension of components of the
bifurcation locus. The proof is as for lemma \ref{fundamental},
where it is not used the fact that the parameters space has dimension 1.
\end{proof}

\begin{lemma}
\label{fundamentalcor2}
Given a bifurcation diagram exhibiting a bifurcation line
$\cal B$, 
there exists a family of gradient vector
fields in a neighbourhood of $\cal B$ having only two saddles and with 
$\cal B$ as associated bifurcation diagram.
\end{lemma}
\begin{proof}
By lemma \ref{crpgrf} we choose a family of functions having two saddles 
points and we apply corollary \ref{fundamentalcor}.
\end{proof}

The following corollary makes global the result of lemma
\ref{fundamentalcor2}.

\begin{corollary}
\label{fundamentalcor3}
Given an allowed bifurcation diagram ${\cal B}$, there exists a family of
gradient vector fields having ${\cal B}$ as
associated bifurcation diagram.
\end{corollary}
\begin{proof}
It is enough to apply lemma \ref{fundamentalcor2} in a neighbourhood
of each component ${\cal B}_{ij}$ of the bifurcation locus.
\end{proof}

The second step is to prove that the family of vector fields of
corollary \ref{fundamentalcor3} can be chosen depending linearly on
the parameter $x$. 

\begin{definition}
A bifurcation diagram ${\cal M}$ is a subset of a bifurcation diagram
${\cal N}$ if the bifurcation locus of ${\cal M}$ is a subset of
the bifurcation locus of ${\cal N}$.
\end{definition}

\begin{definition}
Two bifurcation diagrams ${\cal M}$ and ${\cal N}$ 
are equivalent if there exists a diffeomorphism
of $\R^2$ mapping caustic and bifurcation locus of ${\cal M}$ onto
those of ${\cal N}$.
\end{definition}

\begin{theorem}
\label{gradientth}
Let $f:\R^2\rightarrow\R$ be the generating function of a 
Lagrangian submanifold $L$, suppose $0\in\R^2$ is a critical point
of $f$ and $W$ is a compact neighbourhood of 0.
Given a bifurcation diagram ${\cal M}$ containing a caustic $K$
and a bifurcation locus ${\cal B}$,
such that the number of connected components of
$f(W)\setminus({\cal B}\cup K)$ is finite, ${\cal M}$ 
is allowed and $K$ is diffeomorphic to the caustic of a small perturbation 
of $f$ in $W$,
then, if $W$ is sufficiently small, 
there exists a generating function $\tilde{f}=f+f'$,
such that $f'$ is supported on $W$,
and whose associated bifurcation diagram, restricted to $W$, contains a
subdiagram equivalent to ${\cal M}$ restrictred to $W$.
\end{theorem}
\begin{proof}
Let $U_i$ be the connected components of $f(W)\setminus({\cal B}\cup K)$:
note that $U_i$ is open and choose
a point $x_i\in U_i$.
The bifurcation diagram ${\cal M}$, being allowed, prescribes  
the classes of orbital equivalence of the phase portrait
of $\nabla\tilde{f}_{x_i}$ in each subset $U_i$, so we define a function
$\tilde{f}_{x_i}$ such that number and nature of its critical points
and behaviour of gradient lines joining each pair of 
these critical points are as 
assigned by 
${\cal M}$. The function $\tilde{f}_{x_i}$ is constructed as $F(~,t)$ in
the proof of proposition \ref{fundamental}, so, for $\tilde{f}_{x_i}$
to satisfy the required conditions, it is enough to define it in 
a neighbourhood
$V_i$ of the chosen critical points and relevant gradient lines.
Observe that we can assume 
$V_i\cap V_j=\emptyset$ 
for $i\neq j$, since the number of $V_i$'s is finite.
Define now $\tilde{f}$ on $\cup U_i$ as
$\tilde{f}(y)=\tilde{f}_{x_i}(y)+x_iy$ and extend it to the whole $W$.
Note that for every $\epsilon>0$, since $\nabla f(0)=0$ and the
conditions $\tilde{f}$ has to satisfy concern its gradient 
$\nabla\tilde{f}$, if $W$ is sufficiently small, then 
$|f'|<\epsilon$, where $|~|$ is a quasi-norm associated with the 
Whitney topology of $C^\infty(\R^2)$. 
Observe also that, since $\nabla\tilde{f}_{x_i}$ is structurally
stable on $V_i$, then there exists a neighbourhood $U^{'}_i$ of $x_i$ in $U_i$
such that $\nabla{\tilde{f}}_{x}$ is orbitally equivalent to 
$\nabla\tilde{f}_{x_i}$ for all $x\in U^{'}_i$.
The function $\tilde{f}$ has
the required properties: indeed, by choosing $\epsilon$ sufficiently small, 
the caustic of $\tilde{f}$ is diffeomorphic to the caustic $K$ in
${\cal M}$; moreover, take a path 
$c:[-1,1]\rightarrow W$ such that
$c(0)=x_i$, $c(1)=x_j$ and $c(t)\in{\cal B}$ for some $t\in[-1,1]$, 
then, as in the proof of proposition \ref{codimbif}, 
there exists a point $t'\in[-1,1]$ such that $c(t')$ is
a bifurcation point for $\tilde{f}$. 
\end{proof}


\begin{thebibliography}{100}\frenchspacing\small

\bibitem{AN} A.A. Andronov, E.A. Leontovich, I.I. Gordon, A.G. Maier,
\emph{Theory of bifurcations of dynamic sustems on a plane,} IPST,
Jerusalem (1971).

\bibitem{AP} D. Arinkin, A. Polishchuk,
\emph{Fukaya category and Fourier transform,}
{\tt math.AG/9811023.}

\bibitem{AGZV} V.I. Arnold, S.M. Gusein-Zade, A.N. Varchenko,
\emph{Singularities of differentiable maps Volume I}, Birkh\"auser,
Boston (1985).
 
\bibitem{BL} Th. Br\"ocker, L. Lander, \emph{Differentiable germs and
  catastrophes,} Cambridge University Press, Cambridge (1975).

\bibitem{BMP1} U. Bruzzo, G. Marelli, F. Pioli, \emph{A Fourier transform for
sheaves on real tori: Part I: the equivalence $Sky(T) \cong
Loc(\hat{T}$,} \rm J. Geo. Phys. {\bf 39} (2001), 174-182.

\bibitem{BMP2} U. Bruzzo, G. Marelli, F. Pioli, \emph{A Fourier transform for
sheaves on real tori: Part II: Relative theory,} 
\rm J. Geo. Phys. {\bf 41} (2002), 312-329.

\bibitem{DM} M. Demazure, \emph{Bifurcations and catastrophes,}
  Springer, Berlin (2000).

\bibitem{F1} K. Fukaya, \emph{Mirror symmetry of Abelian
varieties and multi-theta
functions,} (2000). Available from the web page
{\tt http://www.math.kyoto-u.ac.jp/\~{}fukaya/abelrev.pdf}

\bibitem{F2} K. Fukaya, \emph{Multivalued Morse theory, asymptotics
analysis and mirror symmetry,} (2002). Available from the web page
{\tt http://www.math.kyoto-u.ac.jp/\~{}fukaya/fukayagrapat.dvi} 

\bibitem{GG} M. Golubitsky, V. Guillemin, \emph{Stable mappings and
  their singularities,} Springer-Verlag, New York (1973)

\bibitem{GH} J. Guckenheimer, P. Holmes, \emph{Nonlinear oscillations,
  dynamical systems, and bifurcations of vector fields,} Springer, New
  York (1983).

\bibitem{KH} A. Katok, B. Hasselblatt, \emph{Introduction to the
  theory od dynamical systems,} Cambridge University Press, Cambridge (1995).

\bibitem{LYZ} N.C. Leung, S.-T. Yau, E. Zaslow, \emph{From special
Lagrangian to Hermitian-Yang-Mills via Fourier-Mukai transform,} {\tt
  math.DG/0005118.}

\bibitem{M} G. Marelli, \emph{Two-dimensional Lagrangian singularities 
and bifurcations of gradient lines II,} pre-print.

%
\bibitem{Mi} J. Milnor, \emph{Lectures on the h-cobordism theorem,}
  Princeton University Press, Princeton (1965).

\bibitem{R} C. Robinson, \emph{Dynamical Systems,} CRC Press, Boca
  Raton (1995). 

\bibitem{V} V.A. Vassilyev, \emph{Lagrange and Legendre characteristic
classes,} Gordon and Breach Science Publishers, New York (1988). 

\end{thebibliography}
\end{document}